\documentclass[a4paper,12pt]{article}

\usepackage{amsmath}
\usepackage{a4wide}
\usepackage{amsfonts}
\newenvironment{pf}{\textbf{Proof:}}{\hspace{\stretch{1}}$\square$}

\newtheorem{thm}{Theorem}
\newtheorem{df}{Definition}
\newtheorem{pr}{Proposition}
\newtheorem{lem}{Lemma}
\newtheorem{dthm}{Definition-Theorem}
\newtheorem{clm}{Claim}

\title{Markov Chains Approximation of Jump-Diffusion Quantum Trajectories}
\author{Cl\'ement PELLEGRINI\\
\vspace{-0,3cm}\scriptsize{Institut C.Jordan}\\
\vspace{-0,3cm}\scriptsize{Universit\'e C.Bernard, Lyon 1}\\
\vspace{-0,3cm}\scriptsize{21, av Claude Bernard}\\
\vspace{-0,3cm}\scriptsize{69622 Villeurbanne Cedex}\\
\vspace{-0,3cm}\scriptsize{France}\\
\vspace{-0,3cm}\scriptsize{e-mail: pelleg@math.univ-lyon1.fr}}

\begin{document}

\maketitle

\begin{abstract}``Quantum trajectories'' are solutions of stochastic differential
equations also called Belavkin or Stochastic Schr\"odinger
Equations. They describe random phenomena in quantum measurement
theory. Two types of such equations are usually considered, one is
driven by a one-dimensional Brownian motion and the other is
driven by a counting process. In this article, we present a way to
obtain more advanced models which use jump-diffusion stochastic
differential equations. Such models come from solutions of
martingale problems for infinitesimal generators. These generators
are obtained from the limit of generators of classical Markov
chains which describe discrete models of quantum trajectories.
Furthermore, stochastic models of jump-diffusion equations are
physically justified by proving that their solutions can be
obtained as the limit of the discrete trajectories.

\end{abstract}

\section{Introduction}

In quantum mechanics, many recent investigations make a heavy use
of Quantum Trajectory Theory with wide applications in quantum
optic or in quantum information (cf \cite{MR2271425}). A quantum
trajectory is a solution of a stochastic differential equation
which describes the random evolution of quantum systems undergoing
continuous measurement. These equations are called Stochastic
Schr\"odinger Equations or Belavkin Equations (see \cite{SSEMB}).

The result of a measurement in quantum mechanic is inherently
random, as is namely expressed by the axioms of the theory. The
setup is as follows. A quantum system is characterized by a
Hilbert space $\mathcal{H}$ (with finite or infinite dimension)
and an operator $\rho$, self-adjoint, positive, trace class with
$Tr[\rho]=1$. This operator is called a ``state'' or a ``density
matrix''. The measurable quantities (energy, momentum,
position...) are represented by the self-adjoint operators on
$\mathcal{H}$ and are called ``observable'' of the system. The
accessible data are the values of the spectrum of the observable.
In finite dimension for example, if $A=\sum_{i=0}^{p}\lambda_iP_i$
denotes the spectral decomposition of an observable $A$, the
observation of an eigenvalue $\lambda_i$, in the state $\rho$, is
random and it is obtained with probability:
\begin{equation}\label{proba}
P_{\rho}[\textrm{to observe}\, \lambda_i]=Tr[\rho\,\,P_i].
\end{equation}

Besides, conditionally to the result, the reference state of the
system is modified. If we have observed the eigenvalue
$\lambda_i$, then the principle called "Wave Packet Reduction"
imposes the state $\rho$ to collapse to the new reference state
\begin{equation}\label{reduction}
\rho_i^1=\frac{P_i\rho P_i}{Tr[\rho\,P_i]}.
\end{equation}
 Quantum Trajectory
Theory is then the study of the modification of the state of a
system undergoing a sequence of measurements. In this way, with
the fact that $P_iP_j=0$ if $i\neq j$, a second measurement of the
same observable $A$, in the state $\rho_i^1$, should give
$$P_{\rho_i^1}[\textrm{to observe}\, \lambda_i]=1.$$
The principle $(\ref{reduction})$ imposed the new state to be
$\rho^2_i=\rho^1_i$. It means that after one measurement, the
information contained in the system is destroyed in the sense that
the evolution is stopped.

Actually, in physics applications, a model of indirect measurement
is used in order to not destroy the dynamic. The physical setup is
the one of interaction between a  small system (atom) and a
continuous field (environment). By performing a continuous time
quantum measurement on the field, after the interaction, we get a
partial information of the evolution of the small system without
destroying it.

This partial information is governed by stochastic models of
Belavkin equations. In the literature, there are essentially two
different evolutions.
\begin{enumerate}
 \item{If $(\rho_t)$ designs the state of the system, then
one evolution is described by a diffusive equation:
\begin{equation}\label{diffff}
d\rho_t=L(\rho_t)dt+\left[\rho_t C^\star+
C\rho_t-Tr\left[\rho_t(C+C^\star) \right]\rho_t\right]dW_t,
\end{equation}
where $W_t$ describes a one-dimensional Brownian motion.}

\item {The other is given by a stochastic differential equation
driven by a counting process:
\begin{equation}\label{jump-equation}
d\rho_t=L(\rho_t)dt+\left[\frac{\mathcal{J}(\rho_t)}{Tr[\mathcal{J}(\rho_t)]}-\rho_t\right]
(d\tilde{N}_t-Tr[\mathcal{J}(\rho_t)]dt),
\end{equation}
where $\tilde{N}_t$ is a counting process with
intensity $\int_0^t Tr[\mathcal{J}(\rho_s)]ds$.}
\end{enumerate}
Equations $(\ref{diffff})$ and $(\ref{jump-equation})$ are called
classical Belavkin Equations. The solutions of these equations are
called ``continuous quantum trajectories''. Such models describe
essentially the interaction between a two-level atom and a spin
chain (\cite{diffusion},\cite{poisson}). More complicated models
(with high degree of liberty) are given by diffusive evolution
with jump described by jump-diffusion stochastic differential
equations.

Even in the classical cases $(\ref{diffff})$ and
$(\ref{jump-equation})$, Belavkin equations pose tedious problems
in terms of physical and mathematical justifications. First
rigorous results are due to Davies \cite{OQS} which has described
the evolution of a two-level atom undergoing a continuous
measurement. Heuristic rules can be used to obtain classical
Belavkin equations $(\ref{diffff})$ and $(\ref{jump-equation})$. A
rigorous way to obtain these stochastic models is to use Quantum
Filtering Theory (\cite{SSEMB}). Such approach needs high analytic
technologies as Von Neumann algebra and conditional expectation in
operator algebra. The physical justification in this way is far
from being obvious and clear. Furthermore technical difficulties
are increased by introducing more degrees of liberty and such
problems are not really treated.

A more intuitive approach consists in using a discrete model of
interaction called ``Quantum Repeated interactions''. Instead of
considering an interaction with a continuous field, the
environment is represented as an infinite chain of identical and
independent quantum system (with finite degree of liberty). Each
part of the environment interacts with the small system during a time interval
of length $h$. After each interaction, a quantum measurement of an
observable of the field is performed. As regards the small system,
the result of observation is rendered by a random modification of
its reference state in the same fashion of $(\ref{reduction})$.
Then the results of measurements can be described by classical
Markov chains called ``discrete quantum trajectories''. Discrete
quantum trajectories depend on the time interaction $h$. By using
Markov Chain Approximation Theory (using notion of infinitesimal
generators for Markov processes), stochastic models for Continual
Quantum Measurement Theory can be justified as continuous time
limit of discrete trajectories. These models are mathematically justified as follows. Infinitesimal
generators are obtained as limit $(h\rightarrow0)$ of generators
of the Markov chains. These limit generators give then rise to
general problems of martingale (\cite{LTS},\cite{CSPM}). In this
article, we show that such problems of martingale are solved by
solution of particular jump-diffusion stochastic differential
equations, which should model continuous time measurement theory.
This approach and these models are next physically justified by proving that the solutions of
these SDEs can be obtained naturally as a limit (in distribution)
of discrete quantum trajectories.\\

This article is structured as follows.

Section $1$ is devoted to the description of the discrete model of
quantum repeated interactions with measurement. A probability space is defined to give
account of the random character and the Markov chain property of discrete quantum trajectories. Next we shall focus on the dependence on $h$ for these Markov chains and we introduce asymptotic assumption in order to come into the question of convergence.

In Section $2$, by using Markov chain approximation technics, we
obtain continuous time stochastic models as limits of discrete
quantum trajectories. We compute natural infinitesimal generators
of Markov chains; these generators also depend on the time
interaction $h$. Therefore we obtain infinitesimal generators as
limit ($h\rightarrow0$) of those. It gives then rise to general
problems of martingale which are solved by jump-diffusion
stochastic differential equations.

Finally in Section $3$, we show that discrete quantum trajectories
converge in distribution to the solution of stochastic
differentials equations described in Section $2$. The stochastic
model of jump-diffusion equations is then physically justified as
the limit of this concrete physical procedure.

\subsection{Discrete Quantum Trajectories}

\subsection{Quantum Repeated Measurements}

This section is devoted to make precise the mathematical model of
indirect measurement and the principle of ``Quantum Repeated
Interactions''. Such model is highly used in physical applications
in quantum optics or in quantum information (see Haroche
\cite{MR2271425}). Let us start by describing the interaction
model without measurement.

 A small system is in contact with an infinite
chain of identical and independent quantum systems. Each copy of
the chain interacts with the small system during a defined time
$h$. A single interaction is described as follows.

The small system is represented by the Hilbert space
$\mathcal{H}_0$ equipped with the state $\rho$. A copy of the
environment is described by a Hilbert space $\mathcal{H}$ with a
reference state $\beta$. The compound system describing the
interaction is given by the tensor product
$\mathcal{H}_{0}\otimes\mathcal{H}$. The evolution during the
interaction is given by a self-adjoint operator $H_{tot}$ on the
tensor product. This operator is called the total Hamiltonian. Its
general form is
$$H_{tot}=H_0\otimes I+I\otimes H+H_{int}$$
where the operators $H_0$ and $H$ are the free Hamiltonian of each
system. The operator $H_{int}$ represents the Hamiltonian of
interaction. This defines the unitary-operator
$$U=e^{ih\,H_{tot}}$$
and the evolution of states of $\mathcal{H}_0\otimes\mathcal{H}$,
in the Schr\"odinger picture, is given by
$$\rho\mapsto U\,\rho\, U^\star.$$
 After this first interaction, a second copy of $\mathcal{H}$ interacts with
$\mathcal{H}_0$ in the same fashion and so on.

  As the chain is supposed to be infinite, the
Hilbert space describing the whole sequence of interactions is
 \begin{equation}
\mathbf{\Gamma}=\mathcal{H}_0\otimes\bigotimes_{k\geq
1}\mathcal{H}_k
\end{equation}
where $\mathcal{H}_k$ denotes the $k$-th copy of $\mathcal{H}$. The countable tensor product $\bigotimes_{k\geq 1}\mathcal{H}_k$
means the following. Consider that $\mathcal{H}$ is of finite
dimension and that $\{X_0,X_1,\ldots,X_n\}$ is a fixed orthonormal
basis of $\mathcal{H}$. The orthogonal projector on
$\mathbb{C}X_0$ is denoted by $\vert X_0\rangle\langle X_0\vert$.
This is the ground state (or vacuum state) of $\mathcal{H}$. The
tensor product is taken with respect to $X_0$ (for details, see
\cite{FRTC}).

\textbf{Remark}: A vector $Y$ in a Hilbert space $\mathcal{H}$ is
represented by the application $\vert Y\rangle$ from $\mathbb{C}$
to $\mathcal{H}$ which acts with the following way $\vert
Y\rangle(\lambda)=\vert \lambda Y\rangle$. The linear form on
$\mathcal{H}$ are represented by the operators $\langle Z\vert$
which acts on the vector $\vert Y\rangle$ by $\langle Z\vert\vert
Y\rangle=\langle Z,Y\rangle$, where $\langle\,,\rangle$ denotes the
scalar product of $\mathcal{H}$.

The unitary evolution describing the k-th interaction is given by the operator
$U_k$ which acts as $U$ on
$\mathcal{H}_0\otimes\mathcal{H}_k$, whereas it acts as the
identity operator on the other copies of $\mathcal{H}$. If $\rho$ is a state on
$\mathbf{\Gamma}$, the effect of the $k$-th interaction is:
$$\rho\mapsto U_k\,\rho\, U_k^\star$$
 Hence the result of the $k$ first interactions is described by the operator $V_k$ on
$\mathcal{B}(\mathbf{\Gamma})$ defined by the recursive formula:
 \begin{equation}\left\{\begin{array}{ccc} V_{k+1}&=&U_{k+1}V_k \\
  V_0&=&I\end{array}\right.\end{equation}
  and the evolution of states is then given, in the Schr\"odinger picture, by:
\begin{equation}\label{II}\rho\mapsto V_k\,\rho\, V_k^\star.
 \end{equation}
We present now the indirect measurement principle. The idea is to perform a measurement of
an observable of the field after each interaction.

 A measurement of an observable of $\mathcal{H}_k$ is modelled as follows. Let $A$ be any
observable on $\mathcal{H}$, with spectral decomposition $A=\sum_{j=0}^{p}\lambda_j P_j$. We consider its natural
ampliation on $\mathbf{\Gamma}$:
\begin{equation}\label{ampliation}
A^k:=\bigotimes_{j=0}^{k-1}I\otimes A\otimes\bigotimes_{j\geq k+1}I.
\end{equation}
 The result of the measurement of $A^k$ is random, the accessible data are its eigenvalues. If $\rho$ denotes the reference state of $\mathbf{\Gamma}$, the observation of $\lambda_j$ is obtained with probability
 $$P[\textrm{to observe}\,\,\lambda_j]=Tr[\,\rho\,P^k_j\,],\,\,\,\,j\in\{0,\ldots,p\},$$
 where $P_j^k$ is the ampliation of $P_j$ in the same way as $(\ref{ampliation})$. If we have observed the eigenvalue $\lambda_j$, the ``wave packet reduction'' imposes that the state after measurement is
 $$\rho_j=\frac{P^k_j\,\rho\,P^k_j }{Tr[\,\rho\,P^k_j\,]}\,.$$

 $\textbf{Remark}$: This corresponds to the new reference state depending on the result of the observation. Another measurement of the observable $A^k$ (with respect to this new state) would give $P[\textrm{to observe}\,\,\lambda_j]=1$ (because $P_iP_j=0$ if $i\neq j$). This means that only one measurement after each interaction gives a significant information. We recover the phenomena expressed in the introduction. This justifies the principle of repeated interactions.\\

 The repeated quantum measurements are the combination of the previous description and the successive interactions $(\ref{II})$. After each interaction, the measurement procedure involves a random modification of the system. It defines namely a sequence of random states which is called ``discrete quantum trajectory''.

 The initial state on $\mathbf{\Gamma}$ is chosen to be
 $$\mu=\rho\otimes\bigotimes_{j\geq 1}\beta_j$$
 where $\rho$ is some state on $\mathcal{H}_0$ and each $\beta_i=\beta$ is a fixed state on $\mathcal{H}$. We denote by $\mu_k$ the new state after $k$ interactions, that
is: $$\mu_k=V_k\,\mu\, V_k^\star.$$

The probability space describing the experience of repeated
measurements is $\Omega^{\mathbb{N}^\star}$, where
$\Omega=\{0,\ldots,p\}$. The integers $i$ correspond to the
indexes of the eigenvalues of $A$. We endow
$\Omega^{\mathbb{N}^\star}$ with the cylinder $\sigma$-algebra
generated by the sets:
 $$\Lambda_{i_1,\ldots,i_k}=\{\omega\in\Omega^{\mathbb{N}^\star}/\omega_1=i_1,\ldots,\omega_k=i_k
\}.$$
 The unitary operator $U_j$ commutes with all $P^k$,
for any $k$ and $j$ with $k<j$. For any set $\{i_1,\ldots,i_k\}$,
we can define the following non normalized state
\begin{eqnarray*}\tilde{\mu}(i_1,\ldots,i_k)&=&(I\otimes P_{i_1}\otimes
\ldots\otimes P_{i_k}\otimes I\ldots) \,\,\mu_k\,\,(I\otimes
P_{i_1}\otimes \ldots\otimes P_{i_k}\otimes
I\ldots)\\
&=&(P^k_{i_k}\ldots P^1_{i_1})\,\,\mu_k\,\,(P^1_{i_1}\ldots P^k_{i_k}).
\end{eqnarray*}
 It is the non-normalized state which corresponds to the successive observation of the eigenvalues $\lambda_{i_1},\ldots,\lambda_{i_k}$ during the $k$ first measurements. The probability to observe these eigenvalues is
 $$P[\Lambda_{i_1,\ldots,i_k}]=P[\textrm{to observe}\,\,(\lambda_{i_1},\ldots,\lambda_{i_k})]=Tr[\tilde{\mu}(i_1,\ldots,i_k)].$$
 This way, we define a probability measure on the cylinder sets of $\Omega^{\mathbb{N}^\star}$ which satisfies the Kolmogorv Consistency Criterion. Hence it defines a unique probability measure on $\Omega^{\mathbb{N}^\star}$. The discrete quantum trajectory on $\mathbf{\Gamma}$ is then given by the following random sequence of states:
 $$\begin{array}{cccc}\tilde{\rho}_k: & \Omega^{\mathbb{N}^\star} &\longrightarrow&
\mathcal{B}(\mathbf{\Gamma})\\
 & \omega & \longmapsto & \tilde{\rho}_k(\omega_1,\ldots,
\omega_k)=\frac{\tilde{\mu}(\omega_1,\ldots,
\omega_k)}{Tr[\tilde{\mu}(\omega_1,\ldots, \omega_k)]}
\end{array}$$
This next proposition follows from the construction and the remarks
above.
\begin{pr}
Let $(\tilde{\rho}_k)$ be the above random sequence of states. We
have for all $\omega\in\Omega^{\mathbb{N}^\star}$:
$$\tilde{\rho}_{k+1}(\omega)=\frac{P^{k+1}_{\omega_{k+1}}\,U_{k+1}\,\,\tilde{\rho}_k(\omega)\,\,U_{k+1}^\star\,
P^{k+1}_{\omega_{k+1}}}{Tr\left[\,\tilde{\rho}_k(\omega)\,U_{k+1}^\star\,
P^{k+1}_{\omega_{k+1}}\,U_{k+1}\right]}.$$
\end{pr}

 The following theorem is an easy
consequence of the previous proposition.

\begin{thm}
 The discrete quantum trajectory $(\tilde{\rho}_n)_n$ is a Markov chain, with values on the
set of states
  of $\mathcal{H}_{0}\bigotimes_{i\geq 1} \mathcal{H}_{i}$. It is described as follows:
\begin{eqnarray*}
P\left[\tilde{\rho}_{n+1}=\mu/\tilde{\rho}_n=\theta_n,\ldots,\tilde{\rho}_0=\theta_0\right]
=P\left[\tilde{\rho}_{n+1}=\mu/\tilde{\rho}_n=\theta_n\right]
\end{eqnarray*}If $\tilde{\rho}_n=\theta_n$, then the random state $\tilde{\rho}_{n+1}$ takes one
of the values:
$$ \frac{
P_i^{n+1}\,(U_{n+1}\,\theta_n\,U_{n+1}^\star)\,
P_i^{n+1}}{Tr\left[\,(U_{n+1}\,\theta_n\,U_{n+1}^\star)\,
P_i^{n+1}\,\right]}\,\,\,\,\,\,i=0,\ldots,p$$ with probability
$Tr\left[\,(U_{n+1}\,\theta_n\,U_{n+1}^\star)\,P_i^{n+1}\,\right].$
\end{thm}

In general, one is more interested into the reduced state on the small system $\mathcal{H}_0$ only. This state is given by taking a partial trace on $\mathcal{H}_0$. Let us recall what partial trace is. If $\mathcal{H}$ is any Hilbert space, we denote by $Tr_{\mathcal{H}}[W]$ the trace of a trace-class operator $W$ on $\mathcal{H}$.

\begin{dthm}
Let $\mathcal{H}$ and $\mathcal{K}$ be two Hilbert spaces. If $\alpha$ is a state on a tensor product
$\mathcal{H}\otimes\mathcal{K}$, then there exists a unique state
$\eta$ on $\mathcal{H}$ which is characterized by the property
 $$\,\ Tr_{\mathcal{H}}[\,\eta
\,X\,]=Tr_{\mathcal{H}\otimes\mathcal{K}}[\,\alpha\, (X\otimes
I)\,]$$ for all $X \in \mathcal{B}(\mathcal{H})$. This unique
state $\eta$ is called the partial trace of $\alpha$ on
$\mathcal{H}$ with respect to $\mathcal{K}$.
\end{dthm}

 Let $\alpha$ be a state on $\mathbf{\Gamma}$, we denote by
$\mathbf{E}_0(\alpha)$ the partial trace of $\alpha$ on
$\mathcal{H}_0$ with respect to $\bigotimes_{k\geq
1}\mathcal{H}_k$. We define a random sequence of states on
$\mathcal{H}_0$ as follows. For all $\omega$ in
$\Omega^{\mathbb{N}^\star}$, define the discrete quantum
trajectory on $\mathcal{H}_0$
\begin{equation}\label{TrPart}\rho_n(\omega)=\mathbf{E}_0[\tilde{\rho}_n(\omega)].\end{equation}
An immediate consequence of Theorem $1$ is the following result.

\begin{thm}\label{RaS}
The quantum trajectory $(\rho_n)_n$ defined by formula $(\ref{TrPart})$ is a
Markov chain with values in the set of states on $\mathcal{H}_0$.
If $\rho_n=\chi_n$, then $\rho_{n+1}$ takes one of the values:
$$\mathbf{E}_0\left[\frac{(I\otimes P_i)\,U(\chi_n\otimes\beta)U^\star\,(I\otimes
P_i)}{Tr[\,U(\chi_n\otimes\beta)U^{\star}\,(I\otimes
P_i)]}\right]\,\,\,\,\,\,i=0\ldots p$$ with probability
$Tr\left[U(\chi_n\otimes\beta)U^\star\,(I\otimes P_i)\right]$.
\end{thm}

 $\textbf{Remark}$: Let us stress that $$\frac{(I\otimes
P_i)\,U\,(\chi_n\otimes\beta)\,U^\star\,(I\otimes
P_i)}{Tr[\,U\,(\chi_n\otimes\beta)\,U^{\star}\,(I\otimes P_i)]}$$ is a
state on $\mathcal{H}_0\otimes\mathcal{H}$. In this situation, the
notation $\mathbf{E}_0$
denotes the partial trace on $\mathcal{H}_0$ with respect to $\mathcal{H}$. The infinite tensor product $\mathbf{\Gamma}$ is just needed to have a clear description of the repeated interactions and the probability space $\Omega^{\mathbb{N}^\star}$.

It is worth noticing that this Markov chain $(\rho_k)$ depends on the time
interaction $h$. By putting $h=1/n$, we can define
 for all $t>0$ \begin{equation}\label{c}\rho_n(t)=\rho_{[nt]}.\end{equation}
It defines then a sequence of processes $(\rho_n(t))$ and we aim to show next that this sequence of processes converges in distribution $(n\rightarrow\infty)$. As announced in the introduction, such convergence is obtained from the convergence of Markov generators of Markov chains. The following section is then devoted to present these generators for quantum trajectories.

\subsection{Infinitesimal Generators}

In all this section we fix a integer $n$. Let $A$ be an observable
and let $\rho_n(t)$ be the process defined from the quantum
trajectory describing the successive measurements of $A$. In this
section, we investigate the explicit computation of the
Markov generator $\mathcal{A}_n$ of the process
$(\rho_n(t))$ (we will make no distinction between the
infinitesimal generators of the Markov chains $(\rho_k)$ and the
process $(\rho_n(t))$ generated by this Markov chain). For
instance, let us introduce some notation.

Let work with $\mathcal{H}_0=\mathbb{C}^{K+1}$. The set of
operators on $\mathcal{H}_0$ can be identified with $\mathbb{R}^P$
for some $P$ (we have $P=2^{(K+1)^2}$, we will see later that we
do not need to give any particular identification). We set
$E=\mathbb{R}^P$ and the set of states becomes then a compact
subset of $\mathbb{R}^P$ (a state is an operator positive with trace $1$).
We denote by $\mathcal{S}$ the set of states and $E=\mathbb{R}^P$.
 For any state $\rho\in\mathcal{S}$, we define
\begin{eqnarray}\label{Transition}
 \mathcal{L}_i^{(n)}(\rho)&=&\mathbf{E}_0\left[\frac{(I\otimes P_i)\,U(n)(\rho\otimes\beta)U^\star(n)\,(I\otimes
P_i)}{Tr[\,U(n)(\rho\otimes\beta)U^{\star}(n)\,(I\otimes P_i)]}\right]
\,\,i=0\ldots p\nonumber\\
p^i(\rho)&=&Tr[U(n)(\rho\otimes\beta)U^\star(n)I\otimes P_i]
\end{eqnarray}
The operators $\mathcal{L}_i^{(n)}(\rho)$ represent transition
states of Markov chains described in Theorem $(\ref{RaS})$ and the
numbers $p^i(\rho)$ are the associated probabilities.
Markov generators for $(\rho_n(t))$ are then expressed as
follows.

\begin{df}
 Let $(\rho_k)$ be a discrete quantum trajectory obtained from the measurement of an observable $A$ of the form $A=\sum\lambda_i P_i$. Let $(\rho_n(t)$ be the process obtained from $(\rho_k)$ by the expression $(\ref{c})$. Let define $P^{(n)}$ the probability measure which satisfies
\begin{eqnarray}
 &&P^{(n)}[\rho_n(0)=\rho]=1\\
&&P^{(n)}[\rho_n(s)=\rho_k,\,k/n\leq s<(k+1)/n]=1\\
&&P^{(n)}[\rho_{k+1}\in{\mathbf{\Gamma}}\big{/}\mathcal{M}_k^{(n)}]=\Pi_{n}(\rho_k,\mathbf{\Gamma})
\end{eqnarray}
where $\Pi_{n}(\rho,.)$ is the transition function of the Markov
chain $(\rho_k)$ given by
\begin{equation}
 \Pi_n(\rho,\mathbf{\Gamma})=\sum_{i=0}^pp^i(\rho)\delta_{\mathcal{L}^{(n)}_i(\rho)}(\mathbf{\Gamma})
\end{equation}
for all Borel subset $\Gamma\in\mathcal{B}(\mathbb{R}^P)$.

 For all state $\rho\in\mathcal{S}$ and all functions $f\in C^2_c(E)$ (i.e $C^2$ with compact support), we define
\begin{eqnarray}\label{generator}
\mathcal{A}_nf(\rho)&=&n\int(f(\mu)-f(\rho))\Pi_n(\rho,d\mu)\nonumber\\
&=&n\sum_{i=0}^p\big{(}f(\mathcal{L}^{(n)}_i(\rho))-f(\rho)\big{)}p^i(\rho).
\end{eqnarray}
The operator $\mathcal{A}_n$ is called the ``Markov generator'' of the Markov chain $(\rho_k)$ (or for the process $(\rho_n(t))$).
\end{df}

The complete description of the generator $\mathcal{A}_n$ needs
the explicit expression of $\mathcal{L}^{(n)}_i(\rho)$ for all
$\rho$ and all $i\in\{0,\ldots,p\}$. In order to establish this,
we need to compute the partial trace operation $\mathbf{E}_0$ on
the tensor product $\mathcal{H}_0\otimes\mathcal{H}$. A judicious
choice of basis for the tensor product allow to make computations
easier.

 Let $\mathcal{H}_0=\mathbb{C}^{K+1}$ and let
$(\Omega_0,\ldots,\Omega_K)$ be any orthonormal basis of
$\mathcal{H}_0$. Recall that $(X_0,\ldots,X_N)$ denotes
 an orthonormal basis of $\mathcal{H}$. For the tensor product we
choose the basis \begin{eqnarray*}\mathcal{B}&=&(\Omega_0\otimes
X_0,\ldots,\Omega_K\otimes X_0,\Omega_0\otimes
X_1,\ldots,\Omega_K\otimes X_1,\ldots,\Omega_0\otimes
X_N,\ldots,\Omega_K\otimes X_N).
\end{eqnarray*} In this basis, any $(N+1)(K+1)\times(N+1)(K+1)$ matrix
$M$ on $\mathcal{H}_0\otimes\mathcal{H}$ can be written by blocks as a $(N+1)\times (N+1)$ matrix $M=\left(M_{ij}\right)_{0\leq i,j\leq N}$ where $M_{ij}$ are operators on $\mathcal{H}_0$. Furthermore we have the following result which allows to compute easily the partial trace.
\begin{clm}
 Let $W$ be a state acting on $\mathcal{H}_0\otimes\mathcal{H}$. If $W=(W_{ij})_{0\leq i,j\leq N}$, is the expression of $W$ in the basis $\mathcal{B}$, where
  the coefficients $W_{ij}$ are operators on $\mathcal{H}_0$,
   then the partial trace with respect to $\mathcal{H}$ is given by the formula:
 $$\mathbf{E}_0[W]=\sum_{i=0}^NW_{ii}.$$
\end{clm}

From this result, we can give the expression of the operators $\mathcal{L}^{(n)}_i(\rho)$.
The reference state of $\mathcal{H}$ is chosen to be the orthogonal projector on $\mathbb{C}X_0$, that is,
 with physical notations
 $$\beta=\vert X_0\rangle\langle X_0\vert.$$
 This state is called the
ground state (or vacuum state) in quantum physics. From general result of G.N.S representation in $C^\star$ algebra, it is worth noticing that it is not a restriction. Indeed such representation allows to identify any quantum system $(\mathcal{H},\beta)$ with another system of the form $(\mathcal{K},\vert X_0\rangle\langle X_0\vert)$ where $X_0$ is the first vector of an orthonormal basis of a particular Hilbert space $\mathcal{K}$ (see \cite{MR1468229} for details).

 The unitary operator $U(n)$ is described by blocks as
  $U(n)=\left(U_{ij}(n)\right)_{0\leq i,j\leq N}$ where the coefficients $U_{ij}$ are $(K+1)\times (K+1)$ matrices
   acting on $\mathcal{H}_0$. For $i\in\{0,\ldots,p\}$, we denote $P_i=(p_{kl}^i)_{0\leq k,l\leq N}$ the eigen-projectors of the observable $A$. Hence the non-normalized states $\mathbf{E}_0[I\otimes
P_i\,U(h)(\rho\otimes\beta)U(h)^\star\,I\otimes P_i]$ and the probabilities $p^i(\rho)$ satisfy
\begin{eqnarray}\label{expr}
\mathbf{E}_0[I\otimes P_i\,U(n)(\rho\otimes\beta)U(n)^\star\,I\otimes
P_i]&=&\sum_{0\leq k,l\leq N}\,p_{kl}^i\,U_{k0}(n)\rho U_{l0}^\star(n)\nonumber\\
p^i(\rho)&=&\sum_{0\leq k,l\leq N}\,p_{kl}^i Tr\left[U_{k0}(n)\rho U_{l0}^\star(n)\right].
\end{eqnarray}
By observing that the operator $\mathcal{L}^{(n)}_i(\rho)$
satisfies
\begin{equation}\label{frac}\mathcal{L}^{(n)}_i(\rho)=\frac{\mathbf{E}_0[I\otimes P_i\,U(n)(\rho\otimes\beta)U(n)^\star\,I\otimes
P_i]}{p^i(\rho)},
\end{equation}
for all $i\in\{0,\ldots,p\}$, we have a complete description of
the generator $\mathcal{A}_n$. In order to consider the limit of
$\mathcal{A}_n$, we present asymptotic assumption for the
coefficient $U_{ij}(n)$ in the following section.

\subsection{Asymptotic Assumption}

The choice of asymptotic for $U(n)=(U_{ij}(n))$ are based on the works of Attal-Pautrat in \cite{FRTC}. They have
 namely shown that the operator process defined for all $t>0$ by
$$V_{[nt]}=U_{[nt]}(n)\ldots
U_1(n),$$ which describes the quantum repeated interactions,
weakly converges (in operator theory) to a process $(\tilde{V}_t)$
satisfying a Quantum Langevin equation. Moreover, this convergence is non-trivian, only if the coefficients
$U_{ij}(n)$ obey to certain normalization. When translated in our
context, it express that there exists operators $L_{ij}$ such that
we have for all $(i,j)\in\{0,\ldots,N\}^2$ (recall $N+1$ is the
dimension of $\mathcal{H}$)
\begin{equation}
 \lim_{n\rightarrow \infty}\,n^{\varepsilon_{ij}}\,\left(U_{ij}(n)-\delta_{ij}I \right)=L_{ij}
\end{equation}
where $\varepsilon_{ij}=\frac{1}{2}(\delta_{0i}+\delta_{0j})$. As the expression $(\ref{expr})$ given
 the expression of $\mathcal{L}_i^{(n)}(\rho)$ only involves the first column of $U(n)$, we only keep
 the following asymptotic
\begin{eqnarray*}
 U_{00}(n)&=&I-\frac{1}{n}L_{00}+\circ\left(\frac{1}{n}\right)\\
U_{i0}(h)&=&\frac{1}{\sqrt{n}}L_{i0}+ \circ\left(\frac{1}{\sqrt{n}}\right)\,\,\,\textrm{for}\,\,i>0
\end{eqnarray*}
Another fact which will be important in the computation of limit generators is the following claim.
\begin{clm}\label{CONDITION}
 The unitary condition implies that there exists a self-adjoint operator $H$ such that:
$$L_{00}=-\left(iH+\frac{1}{2}\sum_{i=1}^NL_{i0}L_{i0}^\star\right)$$
Furthermore we have for all $\rho\in\mathcal{S}$:
$$Tr\left[L_{00}\rho+\rho L_{00}^\star+\sum_{1\leq k\leq N}L_{k0}\rho L_{k0}^\star\right]=0$$
because $Tr[U(n)\rho\, U^\star(n)]=1$ for all $n$.
\end{clm}

We can now apply these considerations to give the asymptotic expression of non-normalized states
 and probabilities given by the expression $(\ref{expr})$. For the non-normalized states, we have
\begin{eqnarray}\label{non}
 &&\mathbf{E}_0[I\otimes P_i\,U(n)(\rho\otimes\beta)U(n)^\star\,I\otimes
P_i]\nonumber\\&=&p^i_{00}\,\rho+\frac{1}{\sqrt{n}}\sum_{1\leq
k\leq N} \left(p^i_{k0}L_{k0}\rho+p^i_{0k}\rho
L_{k0}^\star\right)\nonumber\\&&+\frac{1}{n}
\left[p^i_{00}\left(L_{00}\rho+\rho
L_{00}^\star\right)+\sum_{1\leq k,l\leq N}p^i_{kl}L_{k0} \rho
L_{l0}^\star\right]+\circ\left(\frac{1}{n}\right)
\end{eqnarray}
with probabilities
\begin{eqnarray}\label{pro}
p^i(\rho)&=&Tr\left[I\otimes P_i\,U(h)(\rho\otimes\beta)U(h)^\star\,I\otimes
P_i\right]\nonumber\\&=&Tr\Big[\mathbf{E}_0[I\otimes P_i\,U(h)(\rho\otimes\beta)U(h)^\star\,I\otimes
P_i]\Big]\nonumber\\
&=&p_{00}^i+\frac{1}{\sqrt{n}}Tr\left[\sum_{1\leq k\leq N}
\left(p^i_{k0}L_{k0}\rho+p^i_{0k}\rho L_{k0}^\star\right)\right]\nonumber\\&&+\frac{1}{n}Tr\left[
\left(p^i_{00}\left(L_{00}\rho+\rho L_{00}^\star\right)+\sum_{1\leq k,l\leq N}p^i_{kl}L_{k0}
\rho L_{l0}^\star\right)\right]+\circ\left(\frac{1}{n}\right).
\end{eqnarray}
 The asymptotic expression of $\mathcal{L}^{(n)}_i(\rho)$ given by the expression $(\ref{frac})$
  follows then from $(\ref{non})$ and $(\ref{pro})$.
   Following the fact that $p_{00}^i$ is equal to zero or not, we consider three cases.

\begin{enumerate}

\item{If $p^i_{00}=0$, then we have
\begin{equation}
 \mathcal{L}^{(n)}_i(\rho)=\frac{\sum_{1\leq k,l\leq N}p^i_{kl}L_{k0}\rho L_{l0}^\star+\circ(1)}{Tr[\sum_{1\leq k,l\leq N}p^i_{kl}L_{k0}\rho L_{l0}^\star+\circ(1)]}\mathbf{1}_{Tr\left[\sum_{1\leq k,l\leq N}p^i_{kl}L_{k0}\rho L_{l0}^\star\right]\neq0}+\circ\left(1\right)
\end{equation}}

\item{If $p^i_{00}=1$, then we have
\begin{eqnarray}
 \mathcal{L}^{(n)}_i(\rho)&=&\rho+\frac{1}{n}\left[\left(L_{00}\rho+\rho L_{00}^\star\right)+\sum_{1\leq k,l\leq N}p^i_{kl}L_{k0}\rho L_{l0}^\star\right]\nonumber\\&&-\frac{1}{n}Tr\left[\left(L_{00}\rho+\rho L_{00}^\star\right)+\sum_{1\leq k,l\leq N}p^i_{kl}L_{k0}\rho L_{l0}^\star\right]\rho+\circ\left(\frac{1}{n}\right)
\end{eqnarray}}

\item{If $p^i_{00}\notin\{0,1\}$, then we have
\begin{eqnarray}
&&\mathcal{L}^{(n)}_i(\rho)\nonumber\\&=&\rho+\frac{1}{\sqrt{n}}\Bigg{[}\frac{1}{p^i_{00}}\sum_{1\leq k\leq N}\left(p^i_{k0}L_{k0}
\rho+p^i_{0k}\rho L_{k0}^\star\right)\nonumber\\&&\hspace{1,8cm}-\frac{1}{p^i_{00}}Tr\left(\sum_{1\leq k\leq N}
\left(p^i_{k0}L_{k0}\rho+p^i_{0k}\rho L_{k0}^\star\right)\right)\times\rho\Bigg{]}\nonumber\\&&+\frac{1}{n}\Bigg{[}
\frac{1}{p_{00}^i}\left(p^i_{00}\left(L_{00}\rho+\rho L_{00}^\star\right)+\sum_{1\leq k,l\leq N}p^i_{kl}L_{k0}\rho
 L_{l0}^\star\right)\nonumber\\
 &&\hspace{1cm}+\frac{1}{(p_{00}^i)^2}Tr\left(\sum_{1\leq k\leq N} \left(p^i_{k0}L_{k0}\rho+p^i_{0k}\rho
L_{k0}^\star\right)\right)^2\times\rho\nonumber\\&&\hspace{1cm}-\frac{1}{p_{00}^i}Tr\left(p^i_{00}\left(L_{00}\rho+\rho
L_{00}^\star\right)
 +\sum_{1\leq k,l\leq N}p^i_{kl}L_{k0}\rho L_{l0}^\star\right)\times\rho\nonumber\\
&&\hspace{1cm}-\frac{1}{(p_{00}^i)^2} \sum_{1\leq k\leq N}
\left(p^i_{k0}L_{k0}\rho+p^i_{0k}\rho L_{k0}^\star\right)\times
Tr\left(\sum_{1\leq k\leq N} \left(p^i_{k0}L_{k0}\rho+p^i_{0k}\rho
L_{k0}^\star\right)\right)\Bigg{]}\nonumber\\
&&+\circ\left(\frac{1}{n}\right)
\end{eqnarray}}
\end{enumerate}
It is worth noticing that all the $\circ$ are uniform in $\rho$ because we work in the set of
 states $\mathcal{S}$ which is compact.

With this description, we can now compute the generator limit of
$\mathcal{A}_n$ for any quantum trajectory. Next we can establish
continuous time model for quantum measurement. This is the main
subject of the following section.

\section{Jump-Diffusion Models of Quantum Measurement}

In this section, we show that the limit ($n\rightarrow\infty$) of
generators $\mathcal{A}_n$ of discrete quantum trajectories gives
rise to explicit infinitesimal generators. From martingale problem
technics, we interpret these generators as generators of Markov
processes. Besides we show that these processes are solution of
jump-diffusion stochastic differential equations which are a
generalization of the classical Belavkin equations
$(\ref{diffff})$ and $(\ref{jump-equation})$ presented in
Introduction.

Let us make precise the notion of martingale problem in our
framework (see \cite{LTS},\cite{CSPM},\cite{MR2152242} and
\cite{MR838085} for complete references). We still consider the
identification of the set of states as a compact subset of
$E=\mathbb{R}^P$ for some $P$. Let $\Pi$ be a transition kernel on
$E$, let $a(.)=(a_{ij}(.))$ be a measurable mapping on $E$ with
values in the set of positive semi-definite symmetric $P\times P$
matrices and let $b(.)=(b_i(.))$ be a measurable function from $E$
to $E$. Let $f$ be any $C^2_c(E)$ and let $\rho\in E$. In this
article, we consider infinitesimal generators $\mathcal{A}$ of the
form
\begin{eqnarray}\label{generatoR}
 \mathcal{A}f(\rho)&=&\sum_{i=1}^P b_i(\rho)\frac{\partial f(\rho)}{\partial \rho_i}+\frac{1}{2}\sum_{i,j=1}^Pa_{ij}(\rho)\frac{\partial f(\rho)}{\partial \rho_i\partial \rho_j}\nonumber\\&&+\int_{E} \left[f(\rho+\mu)-f(\rho)-\sum_{i=1}^P\mu_i\frac{\partial f(\rho)}{\partial \rho_i}\right]\Pi(\rho,d\mu)
\end{eqnarray}
The notion of problem of martingale associated with such
generators is expressed as follows.
\begin{df}\label{PM}
 Let $\rho_0\in E$. We say that a measurable stochastic process $(\rho_t)$ on some probability space $(\Omega,\mathcal{F},P)$ is a solution of the martingale problem for $(\mathcal{A},\rho_0)$, if for all $f\in C^2_c(E)$,
\begin{equation}
 \mathcal{M}_t^f=f(\rho_t)-f(\rho_0)-\int_0^t\mathcal{A}f(\rho_{s})ds,\,\,\,t\geq0
\end{equation}
is a martingale with respect to
$\mathcal{F}_t^\rho=\sigma(\rho_s,\,\,s\leq t)$.
\end{df}
It is worth noticing that we must also define a probability space
$(\Omega,\mathcal{F},P)$ to make explicit a solution of a problem
of martingale.

In the following section, we show that Markov generators of discrete quantum trajectory converges to infinitesimal generators of the form $(\ref{generatoR})$.

\subsection{Limit Infinitesimal Generators}

Before to express the proposition which gives the limit infinitesimal generators of $\mathcal{A}_n$
 defined in Section $1$, we define some functions which appears in the limit. For all $i$ and all state
  $\rho\in\mathcal{S}$, set
\begin{eqnarray}
 g_i(\rho)&=&\Bigg{(}\frac{\sum_{1\leq k,l\leq N}p^i_{kl}L_{k0}\rho L_{l0}^\star}{Tr[\sum_{1\leq k,l\leq N}p^i_{kl}L_{k0}\rho L_{l0}^\star]}-\rho\Bigg{)}\nonumber\\
v_i(\rho)&=&Tr\left[\sum_{1\leq k,l\leq N}p^i_{kl}L_{k0}\rho L_{l0}^\star\right]\nonumber\\
 h_i(\rho)&=&\frac{1}{\sqrt{p_{00}^i}}\left[\sum_{1\leq k\leq N}\left(p^i_{k0}L_{k0}
\rho+p^i_{0k}\rho L_{k0}^\star\right)-Tr\Big{[}\sum_{1\leq k\leq N}
\left(p^i_{k0}L_{k0}\rho+p^i_{0k}\rho L_{k0}^\star\right)\Big{]}\rho\right]\nonumber\\
L(\rho)&=&L_{00}\rho+\rho L_{00}^\star+\sum_{1\leq k\leq N}L_{k0}\rho L_{k0}^\star.
\end{eqnarray}

This next proposition concerning limit generators follows from results of asymptotic described in Section $1$.
\begin{pr}\label{defgen}Let $A$ be an observable with spectral decomposition
$A=\sum_{i=0}^p\lambda_iP_i$ where $P_i=(p^i_{kl})_{0\leq k,l\leq
N}$ are its eigen-projectors. Up to permutation of
eigen-projectors, we can suppose that $p_{00}^0\neq0$. We define
the sets
\begin{eqnarray*}I&=&\{i\in\{1,\ldots,p\}/p_{00}^i=0\}\,\,
\textrm{and}\\ J&=&\{1,\ldots,p\}\setminus I.
 \end{eqnarray*}
 Let $(\rho_n^{\mbox{\tiny{J}}}(t))$ be the corresponding quantum trajectory obtained from the measurement of $A$
  and let $\mathcal{A}_n^J$ be its infinitesimal generator (cf Definition $1$). Let $\mathcal{A}^J$ be the limit
  generator (if it exists) of $\mathcal{A}_n^J$. It is described as follows.
\begin{enumerate}
\item{If $I=\{1,\ldots,p\}$, then $p_{00}^0=1$ and $J=\emptyset$, we have for all $f\in C^2_c(E)$:
\begin{equation}
 \lim_{n\rightarrow\infty}\sup_{\rho\in\mathcal{S}}\big{\vert}\mathcal{A}^J_nf(\rho)-\mathcal{A}^Jf(\rho)\big{\vert}
\end{equation}
where $\mathcal{A}^J$ satisfies
\begin{eqnarray}
\mathcal{A}^Jf(\rho)&=&D_{\rho}f\bigg{(}L(\rho)\bigg{)}
+\int_E\bigg{[}f\bigg{(}\rho+\mu\bigg{)}-f(\rho)-D_\rho
f(\mu)\bigg{]}\Pi(\rho,d\mu),
\end{eqnarray}
the transition kernel $\Pi$ being defined as
$$\Pi(\rho,d\mu)=\sum_{i=1}^pv_i(\rho)\delta_{g_i(\rho)}(d\mu).$$}
\item{If $I\neq\{1,\ldots,p\}$, then $p_{00}^0\neq1$ and $J\neq\emptyset$, we have for all $f\in C^2_c(E)$:
\begin{equation}
 \lim_{n\rightarrow\infty}\sup_{\rho\in\mathcal{S}}\big{\vert}\mathcal{A}^J_nf(\rho)-\mathcal{A}^Jf(\rho)\big{\vert}
\end{equation}
where $\mathcal{A}^J$ satisfies
\begin{eqnarray}
 \mathcal{A}^Jf(\rho)&=&D_\rho f(L(\rho))+\frac{1}{2}\sum_{i\in J\bigcup\{0\}}D^2_\rho f(h_i(\rho),h_i(\rho))\nonumber\\&&+\int_E \left[f(\rho+\mu)-f(\rho)-D_\rho f(\mu)\right]\Pi(\rho,d\mu),
 \end{eqnarray}
the transition kernel $\Pi$ being defined as
$$\Pi(\rho,d\mu)=\sum_{i\in I}v_i(\rho)\delta_{g_i(\rho)}(d\mu).$$}
\end{enumerate}
\end{pr}

\begin{pf} Recall that $\mathcal{S}$ is the set of states and it is a compact subset of $E$.
 For any $i\in\{0,\ldots,p\}$ and for any $\rho\in\mathcal{S}$, let compute
$$\lim_{n\rightarrow\infty}\mathcal{A}^J_nf(\rho)$$
For this aim, we use asymptotic results of Section $1$. As was
described, there are three cases.
\begin{enumerate}
\item{Suppose $p_{00}^i=0$, we have,
\begin{eqnarray}
&&\lim_{n\rightarrow\infty}n\left(f(\mathcal{L}^{(n)}_i(\rho))-f(\rho)\right)p^i(\rho)\nonumber\\
&=&\bigg{[}f\bigg{(}\frac{\sum_{1\leq k,l\leq N}p^i_{kl}L_{k0}\rho L_{l0}^\star}{Tr[\sum_{1\leq k,l\leq N}p^i_{kl}L_{k0}\rho L_{l0}^\star]}\bigg{)}-f(\rho)\bigg{]}Tr\left[\sum_{1\leq k,l\leq N}p^i_{kl}L_{k0}\rho L_{l0}^\star\right]
\end{eqnarray}
Moreover, since we have $f\in C^2_c$ and since $\mathcal{S}$ is
compact, the function defined on $\mathcal{S}$ by
$$\bigg{[}f\bigg{(}\frac{\sum_{1\leq k,l\leq N}p^i_{kl}L_{k0}\rho L_{l0}^\star}{Tr[\sum_{1\leq k,l\leq N}p^i_{kl}L_{k0}\rho L_{l0}^\star]}\bigg{)}-f(\rho)\bigg{]}Tr\left[\sum_{1\leq k,l\leq N}p^i_{kl}L_{k0}\rho L_{l0}^\star\right]$$
is uniformly continuous. As a consequence, the asymptotic concerning this case (and the fact that all the $\circ$ are uniform on $\mathcal{S}$ cf Section $1$) implies the uniform convergence.}
\item{Suppose $p_{00}^i=1$, by using the Taylor formula of order one, we have
\begin{eqnarray}
&&\lim_{n\rightarrow\infty}n\left(f(\mathcal{L}^{(n)}_i(\rho))-f(\rho)\right)p^i(\rho)\nonumber\\
&=&D_{\rho}f\bigg{(}\Big{[}\left(L_{00}\rho+\rho L_{00}^\star\right)+\sum_{1\leq k,l\leq N}p^i_{kl}L_{k0}\rho L_{l0}^\star\Big{]}\nonumber\\&&\hspace{1,2cm}-Tr\Big{[}\left(L_{00}\rho+\rho L_{00}^\star\right)+\sum_{1\leq k,l\leq N}p^i_{kl}L_{k0}\rho L_{l0}^\star\Big{]}\rho\bigg{)}
\end{eqnarray}
To obtain the uniform result, we use asymptotic of Section $1$ and
the uniform continuity of $Df$ on $\mathcal{S}$.}
\item{Suppose $P_{00}^i\notin\{0,1\}$. By applying the Taylor formula of order two, we get the convergence
\begin{eqnarray}
&&\sum_{i/p_{00}^i\notin\{0,1\}}\lim_{n\rightarrow\infty}n\left(f(\mathcal{L}^{(n)}_i(\rho))-f(\rho)\right)p^i(\rho)\nonumber\\
&=&\sum_{i/p_{00}^i\notin\{0,1\}}\Bigg{[}D_{\rho}f\bigg{(}\Big{(}p^i_{00}\left(L_{00}\rho+\rho L_{00}^\star\right)+\sum_{1\leq k,l\leq N}p^i_{kl}L_{k0}\rho
 L_{l0}^\star\Big{)}\nonumber\\&&\hspace{2,8cm}-Tr\Big{[}p^i_{00}\left(L_{00}\rho+\rho
L_{00}^\star\right)
 +\sum_{1\leq k,l\leq N}p^i_{kl}L_{k0}\rho L_{l0}^\star\Big{]}\rho\bigg{)}\nonumber\\
&&+\frac{1}{2p^i_{00}}D^2_{\rho}f\bigg{(}\sum_{1\leq k\leq N}\left(p^i_{k0}L_{k0}
\rho+p^i_{0k}\rho L_{k0}^\star\right)-Tr\Big{[}\sum_{1\leq k\leq N}
\left(p^i_{k0}L_{k0}\rho+p^i_{0k}\rho L_{k0}^\star\right)\Big{]}\rho\nonumber\\&&\hspace{1,2cm},\sum_{1\leq k\leq N}\left(p^i_{k0}L_{k0}
\rho+p^i_{0k}\rho L_{k0}^\star\right)-Tr\Big{[}\sum_{1\leq k\leq N}
\left(p^i_{k0}L_{k0}\rho+p^i_{0k}\rho L_{k0}^\star\right)\Big{]}\rho\bigg{)}\Bigg{]}.\nonumber\\
\end{eqnarray}
Let explain more precisely the last equality. When we use the
Taylor formula for each $i$ such that $p_{00}^i\notin\{0,1\}$, the
term
$$G_i(\rho)=\frac{1}{\sqrt{n}}D_\rho f\Bigg{(}\sum_{1\leq k\leq N}\left(p^i_{k0}L_{k0}
\rho+p^i_{0k}\rho L_{k0}^\star\right)-Tr\left(\sum_{1\leq k\leq N}
\left(p^i_{k0}L_{k0}\rho+p^i_{0k}\rho L_{k0}^\star\right)\right)\rho\Bigg{)}$$
appears, but we have
$$\sum_{i/p_{00}^i\notin\{0,1\}}G_i(\rho)=0$$
since
$\sum_{i/p_{00}^i\notin\{0,1\}}p_{k0}^i=\sum_{i/p_{00}^i\notin\{0,1\}}p_{0k}^i=\sum_{i=0}^p
p_{0k}=\sum_{i=0}^p p_{k0}=0$ for any $k>0$ (indeed we have
$\sum_{i=0}^pP_i=I$). Furthermore this convergence is uniform for
the same arguments as previously.}
\end{enumerate}
These three convergence allow us to obtain the two different cases
of the proposition. The first case of Proposition $2$ follows from
the first two convergences described above, the second case
follows from the first and the third convergences above. Before to
describe this in details, we have to notice that
$$\sum_{i=0}^p\sum_{1\leq k,l\leq N}p^i_{kl}L_{k0}\rho L_{l0}^\star=\sum_{1\leq k\leq N}L_{k0}\rho L_{k0}^\star$$
since we work with eigen-projectors ($\sum_{i=0}^pP_i=Id$). This
fact will be used several times. Moreover, we have
$$Tr[L(\rho)]=Tr\left[L_{00}\rho+\rho L_{00}^\star+\sum_{1\leq k\leq N}L_{k0}\rho L_{k0}^\star\right]=0$$
because of Claim $2$ in Section $1$ concerning the fact that $U$ is a unitary operator.

Using these facts, in case $p^i_{00}=0$, the limit can be written as
$$\int_E\big{[}f(\rho+\mu)-f(\rho)-D_\rho f(\mu)\big{]}v_i(\rho)\delta_{g_i(\rho)}(d\mu)+ D_\rho f(g_i(\rho))v_i(\rho).$$
Besides, we have
$$D_\rho f(g_i(,\rho))v_i(\rho)=D_\rho f\left(\sum_{1\leq k,l\leq N} p^i_{kl}L_{k0}\rho L^\star_{l0}-Tr\left[\sum_{1\leq k,l\leq N} p^i_{kl}L_{k0}\rho L^\star_{l0}\right]\rho\right).$$
Hence it implies the first case of Proposition $2$. For
$I=\{1,\ldots,p\}$, we get indeed
\begin{eqnarray*}
 \mathcal{A}^J_f(\rho)&=&D_\rho f(L(\rho)-Tr[L(\rho)]\rho)+\int_E[f(\rho+\mu)-f(\rho)-D_\rho f(\mu)]\Pi(\rho,d\mu)\\
&=&D_\rho f(L(\rho))+\int_E[f(\rho+\mu)-f(\rho)-D_\rho
f(\mu)]\Pi(\rho,d\mu).
\end{eqnarray*}

A similar reasonment gives the expression of the infinitesimal
generator in the second case where $I\neq\{1,\ldots,p\}$ and the
proposition is proved.
\end{pf}\\

It is worth noticing that generators $\mathcal{A}^J$ are
generators of type $(\ref{generatoR})$, it suffices to expand the
differential terms $D_\rho f$ and $D^2_\rho f$ in terms of partial
derivatives $\frac{\partial f}{\partial \rho_j}$ and
$\frac{\partial^2 f}{\partial \rho_i\partial \rho_j}$.

In the next section, we present continuous time stochastic models
which follows from problems of martingale for the limit
infinitesimal generators $\mathcal{A}^J$.

\subsection{Solutions of Problem of Martingale}

In all this section, we consider an observable $A$ with spectral
decomposition
\begin{equation}\label{decompo}A=\sum_{i\in I}\lambda_iP_i+\sum_{j\in
J\bigcup0}\lambda_jP_j,\end{equation} where $I$ and $J$ are the
subsets of $\{1,\ldots,p\}$ involved in Proposition $2$. Let
$\mathcal{A}^J$ be the associated limit generator and let
$\rho_0$ be a state. In order to solve the problem of martingale
for $(\mathcal{A}^J,\rho_0)$, by Definition $\ref{PM}$, we have to
define a probability space $(\Omega,\mathcal{F},P)$ and a
stochastic process $(\rho_t^{\mbox{\tiny{J}}})$ such that the
process
\begin{equation}\label{proMart}
\mathcal{M}^f_t=f(\rho_t^{\mbox{\tiny$J$}})-f(\rho_0)-\int_0^t\mathcal{A}^Jf(\rho_s^{\mbox{\tiny$J$}})ds
\end{equation}
is a martingale for the natural filtration of $(\rho_t^{\mbox{\tiny$J$}})$.

A classical way to solve the problem of martingale is to define
the solution through a stochastic differential equation
(\cite{MR2152242},\cite{MR1183654}).

Let define a suitable probability space which satisfies the
martingale problem. Consider $(\Omega,\mathcal{F},P)$ a
probability space which supports a $(p+1)$-dimensional Brownian
motion $W=(W_0,\ldots,W_p)$ and $p$ independent Poisson point
processes $(N_i)_{1\leq i\leq p}$ on $\mathbb{R}^2$ and
independent of the Brownian motion.

 As there are two types of limit generators in Proposition \ref{defgen},
  we define two types of stochastic differential equations in the following
  way. Let $\rho_0$ be an initial deterministic state.
\begin{enumerate}
 \item{In case $J=\emptyset$, we define the following stochastic differential
 equation on $(\Omega,\mathcal{F},P)$
\begin{equation}\label{poissonseul}
 \rho^{\mbox{\tiny$J$}}_t=\rho_0+\int L(\rho^{\mbox{\tiny$J$}}_{s-})ds+\sum_{i=1}^p\int_0^t\int_{\mathbb{R}}g_i(\rho^{\mbox{\tiny$J$}}_{s-})\mathbf{1}_{0<x<v_i(\rho^{\mbox{\tiny$J$}}_{s-})}\left[N_i(dx,ds)-dxds\right].
\end{equation}}
\item{In case $J\neq\emptyset$, we define
\begin{eqnarray}\label{avecdiffusif}
 \rho^{\mbox{\tiny$J$}}_t&=&\rho_0+\int L(\rho^{\mbox{\tiny$J$}}_{s-})ds+\sum_{i\in J\bigcup\{0\}}\int_0^th_i(\rho^{\mbox{\tiny$J$}}_{s-})dW_i(s)\nonumber\\&&+\sum_{i\in I}\int_0^t\int_{\mathbb{R}}g_i(\rho^{\mbox{\tiny$J$}}_{s-})\mathbf{1}_{0<x<v_i\left(\rho^{\mbox{\tiny$J$}}_{s-}\right)}\left[N_i(dx,ds)-dxds\right].
\end{eqnarray}}
\end{enumerate}

In this way of writing, these stochastic differential equations
have a meaning only if the process-solution takes values in the
set of states (in general the term $v_i(\rho)$ is not real for all
operator $\rho$). We must modify the expression in order to
consider such equation in a general way for all process which
takes values in operators on $\mathcal{H}_0$. For all $i$, we
define when it has a meaning:
$$\tilde{g}_i(\rho)=\frac{\sum_{1\leq k,l\leq N}p^i_{kl}L_{k0}\rho L_{l0}^\star}{Re(v_i(\rho))}.$$
Hence we consider the modified stochastic differential equations
\begin{eqnarray}\label{modified2}
\rho^{\mbox{\tiny$J$}}_t&=&\rho_0+\int L(\rho^{\mbox{\tiny$J$}}_{s-})ds+\sum_{i=1}^p\int_0^t\int_{\mathbb{R}}\tilde{g}_i(\rho^{\mbox{\tiny$J$}}_{s-})\mathbf{1}_{0<x<Re(v_i(\rho^{\mbox{\tiny$J$}}_{s-}))}\left[N_i(dx,ds)-dxds\right]
\end{eqnarray}
and
\begin{eqnarray}\label{modified1}
 \rho^{\mbox{\tiny$J$}}_t&=&\rho_0+\int L(\rho^{\mbox{\tiny$J$}}_{s-})ds+\sum_{i\in J\bigcup\{0\}}\int_0^th_i(\rho^{\mbox{\tiny$J$}}_{s-})dW_i(s)\nonumber\\&&+\sum_{i\in I}\int_0^t\int_{\mathbb{R}}\tilde{g}_i(\rho^{\mbox{\tiny$J$}}_{s-})\mathbf{1}_{0<x<Re(v_i(\rho^{\mbox{\tiny$J$}}_{s-}))}\left[N_i(dx,ds)-dxds\right],
\end{eqnarray}
Let $\rho$ be a state. The fact that $Re(v_i(\rho))=v_i(\rho)$ and
$\tilde{g}_i(\rho)=g_i(\rho)$ implies that a solution
$(\rho^{\mbox{\tiny$J$}}_t)$ of the equation $(\ref{modified2})$
(resp $(\ref{modified1})$) is a solution of the equation
$(\ref{poissonseul})$ (resp $(\ref{avecdiffusif})$) when the
process $(\rho^{\mbox{\tiny$J$}}_t)$ takes values in the set of
states.

We proceed in the following way to solve the problem of martingale
$(\ref{proMart})$. Firstly we show that the modified equations
$(\ref{modified2})$ and $(\ref{modified1})$ admit a unique
solution (we will see below that it needs another modification),
secondly we show that solutions of $(\ref{modified2})$ and
$(\ref{modified1})$ can be obtained as limit (in distribution) of
discrete quantum trajectories (cf Section $3$). Finally we show
that the property of being a process valued in the set of states
follows from convergence (cf Section $3$) and we conclude that
solutions of $(\ref{modified2})$ and $(\ref{modified1})$ takes
values in the set of states. Moreover, we show that they are solutions of problem
of martingale $(\ref{proMart})$.

The fact that if solutions of $(\ref{modified2})$ and
$(\ref{modified1})$ takes values in the set of states, they are
solutions of martingale problem $(\ref{proMart})$ is expressed in
the following proposition.

\begin{pr}\label{sol} Let $\rho_0$ be any initial state.

If the modified stochastic differential equations
$(\ref{modified2})$ admits a solution $(\rho^{\mbox{\tiny$J$}}_t)$
which takes values in the set of states, then it is a solution of
the problem of martingale $(\mathcal{A}^J,\rho_0)$ in the case
$I=\{1,\ldots,p\}$.

If the modified stochastic differential equations
$(\ref{modified1})$ admits a solution $(\rho^{\mbox{\tiny$J$}}_t)$
which takes values in the set of states, then it is a solution of
the problem of martingale $(\mathcal{A}^J,\rho_0)$ in the case
$J\neq\emptyset$.

As a consequence, if $\tilde{\mathcal{A}}^J$ designs the
infinitesimal generator of a solution of $(\ref{modified1})$ or
$(\ref{modified2})$, then we have
$\tilde{\mathcal{A}}^Jf(\rho)=\mathcal{A}^Jf(\rho)$ for all state
$\rho$ and all functions $f\in C^2_c$.
\end{pr}

\begin{pf}
Recall we assume that processes take values in the set of
states. For any state $\rho$, we have $Re(v_i(\rho))=v_i(\rho)$
and $\tilde{g}_i(\rho)=g_i(\rho)$ and the part concerning the
generators follows.
 Concerning the martingale problem, it is a consequence of the It\^o formula.
 Let $\rho^{\mbox{\tiny$J$}}_t=(\rho^{\mbox{\tiny$J$}}_1(t),\ldots,\rho^{\mbox{\tiny$J$}}_{P}(t))$ denote the coordinates of
 a solution of $(\ref{poissonseul})$ or $(\ref{avecdiffusif})$ (with identification between the set of operators on $\mathcal{H}_0$ and $\mathbb{R}^P$), we have for all $f\in C^2_c$
\begin{eqnarray}
f(\rho^{\mbox{\tiny$J$}}_t)-f(\rho_0)&=&\sum_{i=1}^P\int_{0}^t\frac{\partial f(\rho^{\mbox{\tiny$J$}}_{s-})}{\partial \rho^{\mbox{\tiny$J$}}_i}d\rho^{\mbox{\tiny$J$}}_i(s)\nonumber\\
&&+\frac{1}{2}\sum_{i,j=1}^P\int_{0}^t\frac{\partial f(\rho^{\mbox{\tiny$J$}}_{s-})}{\partial\rho^{\mbox{\tiny$J$}}_i\partial\rho^{\mbox{\tiny$J$}}_j}d[\rho^{\mbox{\tiny$J$}}_i(s),\rho^{\mbox{\tiny$J$}}_j(s)]^c\nonumber\\
&&+\sum_{0\leq s\leq t}\left[f(\rho^{\mbox{\tiny$J$}}_s)-f(\rho^{\mbox{\tiny$J$}}_{s-})-\sum_{i=1}^P\frac{\partial f(\rho^{\mbox{\tiny$J$}}_{s-})}{\partial\rho^{\mbox{\tiny$J$}}_i}\Delta\rho^{\mbox{\tiny$J$}}_i(s)\right]
\end{eqnarray}
where $[\rho^{\mbox{\tiny$J$}}_i(.),\rho^{\mbox{\tiny$J$}}_j(.)]^c$ denotes the continuous part of $[\rho^{\mbox{\tiny$J$}}_i(.),\rho^{\mbox{\tiny$J$}}_j(.)]$.

Let us deal with the case where $J\neq\emptyset$. If $(e_i)_{1\leq
i\leq P}$ designs the canonical basis of $\mathbb{R}^P$, then we
have $\rho^{\mbox{\tiny$J$}}_i(t)=\langle
\rho^{\mbox{\tiny$J$}}_t,e_i\rangle$ for all $t\neq0$. Hence we
have $d\rho^{\mbox{\tiny$J$}}_i(t)=\langle
d\rho^{\mbox{\tiny$J$}}_t,e_i\rangle$. As a consequence we have
for all $i\in\{1,\ldots,P\}$
\begin{eqnarray}
 \rho^{\mbox{\tiny$J$}}_i(t)&=&\rho_0+\int \langle L(\rho^{\mbox{\tiny$J$}}_{s-}),e_i\rangle ds+\sum_{k\in J\bigcup\{0\}}\int_0^t\left\langle h_k(\rho^{\mbox{\tiny$J$}}_{s-}),e_i\right\rangle dW_k(s)\nonumber\\&&+\sum_{k\in I}\int_0^t\int_{\mathbb{R}}\langle g_k(\rho^{\mbox{\tiny$J$}}_{s-}),e_i\rangle\mathbf{1}_{0<x<v_k(\rho^{\mbox{\tiny$J$}}_{s-})}\left[N_k(dx,ds)-dxds\right].
\end{eqnarray}
It implies that
$$[\rho^{\mbox{\tiny$J$}}_i(t),\rho^{\mbox{\tiny$J$}}_j(t)]^c=\sum_{k\in J\bigcup\{0\}}\int_0^t\left\langle h_k(\rho^{\mbox{\tiny$J$}}_{s-}),e_i\right\rangle\left\langle h_k(\rho^{\mbox{\tiny$J$}}_{s-}),e_i\right\rangle ds$$
since $[W_i(t),W_j(t)]=\delta_{ij}t$. Furthermore, if we set by
$g_k^i(\rho)=\langle g_k(\rho),e_i\rangle$, then we get that the
process
\begin{eqnarray}&&\sum_{0\leq s\leq t}\left[f(\rho^{\mbox{\tiny$J$}}_s)-f(\rho^{\mbox{\tiny$J$}}_{s-})-\sum_{i=1}^P\frac{\partial f(\rho^{\mbox{\tiny$J$}}_{s-})}{\partial\rho_i}\Delta\rho^{\mbox{\tiny$J$}}_i(s)\right]\nonumber\\&&-\sum_{k\in J}\int_0^t\int_{\mathbb{R}}\left[f(\rho^{\mbox{\tiny$J$}}_{s-}+g_k(\rho^{\mbox{\tiny$J$}}_{s-}))-f(\rho^{\mbox{\tiny$J$}}_{s-})-\sum_{i=1}^P\frac{\partial f(\rho^{\mbox{\tiny$J$}}_{s-})}{\partial\rho_i}g_k^i(\rho^{\mbox{\tiny$J$}}_{s-})\right]\mathbf{1}_{0<x\leq v_k(\rho^{\mbox{\tiny$J$}}_{s-})}N_k(dx,ds)\nonumber
\end{eqnarray}
is a martingale. Hence we have
\begin{eqnarray}
 &&\sum_{k\in J}\int_0^t\int_{\mathbb{R}}\left[f(\rho^{\mbox{\tiny$J$}}_{s-}+g_k(\rho^{\mbox{\tiny$J$}}_{s-}))-f(\rho^{\mbox{\tiny$J$}}_{s-})-\sum_{i=1}^P\frac{\partial f(\rho^{\mbox{\tiny$J$}}_{s-})}{\partial\rho^{\mbox{\tiny$J$}}_i}g_k^i(\rho^{\mbox{\tiny$J$}}_s-)\right]\mathbf{1}_{0<x\leq v_k(\rho^{\mbox{\tiny$J$}}_{s-})}N_k(dx,ds)\nonumber\\&&-\sum_{k\in J}\int_0^t\int_{\mathbb{R}}\left[f(\rho^{\mbox{\tiny$J$}}_{s-}+g_k(\rho^{\mbox{\tiny$J$}}_{s-}))-f(\rho^{\mbox{\tiny$J$}}_{s-})-\sum_{i=1}^P\frac{\partial f(\rho^{\mbox{\tiny$J$}}_{s-})}{\partial\rho^{\mbox{\tiny$J$}}_i}g_k^i(\rho^{\mbox{\tiny$J$}}_s-)\right]\mathbf{1}_{0<x\leq v_k(\rho^{\mbox{\tiny$J$}}_{s-})}dxds
\end{eqnarray}
is a martingale because each $N_k$ is a Poisson point process with intensity measure $dx\otimes ds$. Furthermore we have
\begin{eqnarray}
&&\sum_{k\in J}\int_0^t\int_{\mathbb{R}}\left[f(\rho^{\mbox{\tiny$J$}}_{s-}+g_k(\rho^{\mbox{\tiny$J$}}_{s-}))-f(\rho^{\mbox{\tiny$J$}}_{s-})-\sum_{i=1}^P\frac{\partial f(\rho^{\mbox{\tiny$J$}}_{s-})}{\partial\rho^{\mbox{\tiny$J$}}_i}g_k^i(\rho^J_s-)\right]\mathbf{1}_{0<x\leq v_k(\rho^{\mbox{\tiny$J$}}_{s-})}dxds\nonumber\\
&=&\sum_{k\in J}\int_0^t\left[f(\rho^{\mbox{\tiny$J$}}_{s-}+g_k(\rho^{\mbox{\tiny$J$}}_{s-}))-f(\rho^{\mbox{\tiny$J$}}_{s-})-\sum_{i=1}^P\frac{\partial f(\rho^{\mbox{\tiny$J$}}_{s-})}{\partial\rho^{\mbox{\tiny$J$}}_i}g_k^i(\rho^{\mbox{\tiny$J$}}_s-)\right]v_k(\rho^{\mbox{\tiny$J$}}_{s-})ds\nonumber\\
&=&\int_0^t\left[f(\rho^{\mbox{\tiny$J$}}_{s-}+\mu)-f(\rho^{\mbox{\tiny$J$}}_{s-})-
D_{\rho^{\mbox{\tiny$J$}}_{s-}}f(\mu)\right]\Pi(\rho^{\mbox{\tiny$J$}}_{s-},d\mu).
\end{eqnarray}
As the Lebesgue measure of the set of times where
$\rho^J_{s-}\neq\rho^J_s$ is equal to zero, we get that
$$f(\rho^{\mbox{\tiny$J$}}_t)-f(\rho_0)-\int_0^t\mathcal{A}^Jf(\rho^{\mbox{\tiny$J$}}_{s-})ds=f(\rho^{\mbox{\tiny$J$}}_t)-f(\rho_0)-\int_0^t\mathcal{A}^Jf(\rho^{\mbox{\tiny$J$}}_{s})ds$$
and it defines a martingale with respect to the natural filtration of $(\rho_t)$ and the proposition is proved.
\end{pf}

As announced, the first step consists in proving that equations
$(\ref{modified2})$ and $(\ref{modified1})$ admit a unique
solution. Concerning the equation $(\ref{modified2})$, we consider
the following way of writing
\begin{eqnarray}\label{modified22}
\rho^{\mbox{\tiny$J$}}_t&=&\rho_0+\int_0^t L(\rho^{\mbox{\tiny$J$}}_{s-})-\sum_{i=1}^p\tilde{g}_i(\rho^{\mbox{\tiny$J$}}_{s-})Re(v_i(\rho^{\mbox{\tiny$J$}}_{s-}))ds\nonumber\\&&+\sum_{i=1}^p\int_0^t\int_{\mathbb{R}}\tilde{g}_i(\rho^{\mbox{\tiny$J$}}_{s-})\mathbf{1}_{0<x<Re(v_i(\rho^{\mbox{\tiny$J$}}_{s-}))}N_i(dx,ds),
\end{eqnarray}
and in the same way for equation $(\ref{modified1})$ we consider
\begin{eqnarray}\label{modified11}
 \rho^{\mbox{\tiny$J$}}_t&=&\rho_0+\int_0^t L(\rho^{\mbox{\tiny$J$}}_{s-})-\sum_{i\in I}\tilde{g}_i(\rho^{\mbox{\tiny$J$}}_{s-})Re(v_i(\rho^{\mbox{\tiny$J$}}_{s-}))ds+\sum_{i\in J\bigcup\{0\}}\int_0^th_i(\rho^{\mbox{\tiny$J$}}_{s-})dW_i(s)\nonumber\\&&+\sum_{i\in I}\int_0^t\int_{\mathbb{R}}\tilde{g}_i(\rho^{\mbox{\tiny$J$}}_{s-})\mathbf{1}_{0<x<Re(v_i(\rho^{\mbox{\tiny$J$}}_{s-}))}N_i(dx,ds)
\end{eqnarray}

 Sufficient conditions (see \cite{QR}), in order to prove that equations $(\ref{modified22})$ and
  $(\ref{modified11})$ admit a unique solution can be expressed as follows. On the one hand the functions
   $L(.)$, $h_i(.)$ and $\tilde{g}_i(.)Re(v_i(.))$ must be Lipschitz for all $i$. On the other hand
    functions $Re(v_i(.))$ must satisfy that there exists a constant $K$ such that we have, for all $i$ and all operator $\rho$ on $\mathcal{H}_0\simeq\mathbb{R}^P$
\begin{equation}\label{bound}
 \sup_{\rho\in\mathbb{R}^P}\vert Re( v_i(\rho))\vert\leq K
\end{equation}
Actually such conditions (Lipschitz and $(\ref{bound})$) are not
satisfied by the functions $L(.)$, $h_i(.)$, $Re(v_i(.))$ and
$\tilde{g}_i(.)Re(v_i(.))$. However these functions are
$C^\infty$, hence these conditions are in fact locally satisfied.
Therefore a truncature method cam be used to make the
functions $L(.)$, $h_i(.)$ and $\tilde{g}_i(.)Re(v_i(.))$
Lipschitz and functions $Re(v_i(.))$ bounded. It is described
as follows.

Fix $k>0$. A truncature method means that we compose the functions
$L(.)$, $h_i(.)$, $Re(v_i(.)$ and $\tilde{g}_i(.)Re(v_i(.))$ with
a truncature function $\phi^k$ of the form
\begin{eqnarray}
\phi^k(x)&=&(\psi^k(x_i))_{i=1,\ldots, P}\,\,\,\,\textrm{where}\\
\psi^k(x_i)&=&-k\mathbf{1}_{ x_i\leq-k}+x_i\mathbf{1}_{\vert x_i\vert<k}+k\mathbf{1}_{ x_i\neq k}\\
\end{eqnarray}
for all $x=(x_i)\in\mathbb{R}^P$. Hence, if $F$ is any function
defined on $\mathbb{R}^P$, we define the function $F^k$ on
$\mathbb{R}^P$ by
$$F^k(x)=F\left(\phi^k(x)\right)$$
for all $x\in\mathbb{R}^P$. By extension we will note $F^k(\rho)$ when we deal with operators on $\mathcal{H}_0$.

As a consequence, functions $L^k(.)$, $h^k_i(.)$ and $\tilde{g}^k_i(.)Re(v^k_i(.)$ become Lipschitz. Furthermore, as $\phi^k$ is a bounded function, we have
$$\sup_i\sup_{\rho\in\mathbb{R}^P}\vert Re( v_i^k(\rho))\vert\leq K.$$
This theorem follows from these conditions.
\begin{thm}\label{existenceanduniqueness}
 Let $k\in\mathbb{R}^+$ and let $\rho_0$ be any operator on $\mathcal{H}_0$. The following stochastic
 differential equations, in case $J=\emptyset$,
\begin{eqnarray}\label{modified222}
\rho^{\mbox{\tiny$J$}}_t&=&\rho_0+\int_0^t L^k(\rho^{\mbox{\tiny$J$}}_{s-})-\sum_{i=1}^p\tilde{g}^k_i(\rho^{\mbox{\tiny$J$}}_{s-})Re(v^k_i(\rho^{\mbox{\tiny$J$}}_{s-}))ds\nonumber\\&&+\sum_{i=1}^p\int_0^t\int_{\mathbb{R}}\tilde{g}^k_i(\rho^{\mbox{\tiny$J$}}_{s-})\mathbf{1}_{0<x<Re(v_i^k(\rho^{\mbox{\tiny$J$}}_{s-}))}N_i(dx,ds),
\end{eqnarray}
and in case $J\neq\emptyset$
\begin{eqnarray}\label{modified111}
 \rho^{\mbox{\tiny$J$}}_t&=&\rho_0+\int_0^t L^k(\rho^{\mbox{\tiny$J$}}_{s-})-\sum_{i\in I}\tilde{g}^k_i(\rho^{\mbox{\tiny$J$}}_{s-})Re(v^k_i(\rho^{\mbox{\tiny$J$}}_{s-}))ds+\sum_{i\in J\bigcup\{0\}}\int_0^th^k_i(\rho^{\mbox{\tiny$J$}}_{s-})dW_i(s)\nonumber\\&&+\sum_{i\in I}\int_0^t\int_{\mathbb{R}}\tilde{g}^k_i(\rho^{\mbox{\tiny$J$}}_{s-})\mathbf{1}_{0<x<Re(v^k_i(\rho^{\mbox{\tiny$J$}}_{s-}))}N_i(dx,ds)
\end{eqnarray}
 admit a unique solution.

 Let $\overline{\mathcal{A}}_k^J$ be the infinitesimal
generator of the solution of an equation of the form
$(\ref{modified222})$ or $(\ref{modified111})$. For any $f\in
C^2_c$ and any state $\rho$, we have for all $k>1$
$$\overline{\mathcal{A}}_k^Jf(\rho)=\mathcal{A}^Jf(\rho).$$
where $\mathcal{A}^J$ are the infinitesimal generators defined in
Proposition $2$. Furthermore in all cases, the processes defined by
\begin{equation}
 \overline{N}^i_t=\int_0^t\int_{\mathbb{R}}\mathbf{1}_{0<x<Re(v^k_i(\rho^{\mbox{\tiny$J$}}_{s-}))}N_i(dx,ds)
\end{equation}
are counting processes with stochastic intensity
$$t\rightarrow\int_0^t\left[Re(v_i(\rho^{\mbox{\tiny$J$}}_{s-}))\right]^+ds,$$
where $(x)^+=max(0,x)$.
\end{thm}

\begin{pf} The part of this theorem concerning generators is the equivalent
of Proposition $\ref{sol}$. This follows from  Proposition
$\ref{sol}$ and from the fact that, on the set of states
$\mathcal{S}$, we have $\phi^k(\rho)=\rho$ for all $k>1$. Indeed
if $\rho=(\rho_i)_{i=1,\ldots,P}$ is a state, we have
$\vert\rho_i\vert\leq1$ for all $i$.

 The last part of this theorem
follows from properties of Random Poisson Measure Theory and is
treated in details in \cite{poisson} for classical Belavkin
equations $(\ref{jump-equation})$. The proof of Theorem $3$
follows from Lipschitz character and works of Jacod and Protter in
\cite{QR}.

Let us investigate the proof in the case where $J\neq\emptyset$ (the case
$J=\emptyset$ is easy to adapt to this case with a similar proof).

Let us prove that equation $(\ref{modified111})$ admits a unique solution (we suppress the
index $J$ in the solution to lighten the way of writing; we suppress also the index $k$ concerning the truncature). As we have $\,\sup_i\sup_{\rho\in\mathbb{R}^P}\vert Re(v_i(\rho))\vert\leq K$, we can consider Poisson point process on $\mathbb{R}\times[0,K]$; the part concerning the counting process can be then written as
$$\int_0^t\int_{[0,K]}\tilde{g}_i(\rho_{s-})\mathbf{1}_{0<x<Re(v_i(\rho_{s-}))}N(dx,ds)$$
Hence for all $i\in I$ the process
\begin{equation}
\mathcal{N}_i(t)=card\{N_i(.,[0,t]\times[0,K])\}
\end{equation}
defines a classical Poisson process of intensity $K$. As a
consequence, for all $t$, it defines a random sequence
$\left\{(\tau^i_k,\xi^i_k),k\in\{1,\ldots\mathcal{N}_i(t)\}\right\}$
where $\tau^i_k$ designs the jump time of $\mathcal{N}_i(.)$ and
the $\xi^i_k$'s are independent uniform random variables on
$[0,K]$. Consequently, the solution of the stochastic differential
equation is given by
\begin{eqnarray}\label{solution}
 \rho_t&=&\rho_0+\int_0^tL(\rho_{s-})ds-\sum_{i\in I}\int_0^t \tilde{g}_i(\rho_{s-})Re(v_i(\rho_{s-})ds\nonumber\\
&&+\sum_{i\in j\bigcup\{0\}}\int_{0}^th_i(\rho_{s-})dW_i(s)+\sum_{i\in I}\sum_{k=1}^{\mathcal{N}_i(t)}\tilde{g}_i(\rho_{\tau_k^i-})\mathbf{1}_{0<\xi^i_k\leq Re\left(v_i(\rho_{\tau_k^i-}\right)}
\end{eqnarray}
The solution $(\ref{solution})$ is described as follow. Thanks to the Lipschitz property (following from the truncature),
 there exists a unique solution $(\rho_t^1)$ of the equation
\begin{eqnarray}\label{sol1}
 \rho^1_t&=&\rho_0+\int_0^tL(\rho^1_{s-})ds-\sum_{i\in I}\int_0^t \tilde{g}_i(\rho^1_{s-})Re(v_i(\rho^1_{s-})ds\nonumber\\
&&+\sum_{i\in j\bigcup\{0\}}\int_{0}^th_i(\rho^1_{s-})dW_i(s)
\end{eqnarray}
 For all $i\in I$, this solution $(\rho_t^1)$ defines the
function $t\rightarrow Re(v_i(\rho^1_{t-}))$. We define the random
stopping time
$$T_1=\inf\left\{t/\sum_{i\in I}\int_{0}^t\int_{[0,K]}\mathbf{1}_{0<x<Re(v_i(\rho^1_{s-}))}N_i(dx,ds)>1\right\}.$$
By definition of Poisson point processes and by independence, we
have for all $i\neq j$:
$$P\left[\exists\, t\bigg/\int_{0}^t\int_{[0,K]}\mathbf{1}_{0<x<Re(v_i(\rho^1_{s-}))}N_i(dx,ds)=\int_{0}^t\int_{[0,K]}\mathbf{1}_{0<x<Re(v_i(\rho^1_{s-}))}N_j(dx,ds)\right]=0$$
As a consequence at $T_1$, there exists a unique index $i_{T_1}$ such that
$$\int_{0}^{T_1}\int_{[0,K]}\mathbf{1}_{0<x<Re\left(v_{i_{T_1}}(\rho^1_{s-})\right)}N_{i_{T_1}}(dx,ds)=1,$$
and all the other terms concerning the other Poisson point processes (for different indexes of $i_{T_1}$) are equal to zero. Moreover, we have almost surely
$$\int_{0}^{T_1}\int_{[0,K]}\mathbf{1}_{0<x<Re\left(v_{i_{T_1}}(\rho^1_{s-})\right)}N_{i_{T_1}}(dx,ds)=\sum_{k=1}^{\mathcal{N}_{i_{T_1}}}\mathbf{1}_{0<\xi_k^{i_{T_1}}<Re\left(v_{i_{T_1}}(\rho^1_{T_1-})\right)}$$
 We define then the solution of $(\ref{modified222})$ on $[0,T_1]$ in the following way
\begin{equation}\label{We}
 \left\{\begin{array}{cccc}\rho_t&=&\rho^1_t&\textrm{on}\,\,[0,T_1[\\
\rho_{T_1}&=&\tilde{g}_{i_{T_1}}(\rho_{T_1-})
\end{array}\right.
\end{equation}
The operator $\rho_{T_1}$ can then be considered as the initial
condition of the equation $(\ref{sol1})$. Therefore we consider
for $t>T_1$ the process $(\rho^2_t)$ defined by
\begin{eqnarray}\label{sol2}
 \rho^2_t&=&\rho_{T_1}+\int_{T_{1}}^tL(\rho^2_{s-})ds-\sum_{i\in I}\int_{T_{1}}^t \tilde{g}_i(\rho^2_{s-})Re(v_i(\rho^2_{s-})ds\nonumber\\
&&+\sum_{i\in j\bigcup\{0\}}\int_{T_1}^th_i(\rho^2_{s-})dW_i(s).
\end{eqnarray}
In the same fashion as the definition of $T_1$, we can define the random stopping time $T_2$ as
$$T_2=\inf\left\{t>T_1/\sum_{i\in I}\int_{T_1}^t\int_{[0,K]}\mathbf{1}_{0<x<Re(v_i(\rho^2_{s-}))}N_i(dx,ds)>1\right\}.$$
By adapting the expression $(\ref{We})$, we can define the solution
on $[T_1,T_2]$ and so on. By induction, we define then the solution
of $(\ref{modified222})$. The uniqueness comes from the uniqueness
of solution for diffusive equations of type of $(\ref{sol2})$.
Because of the fact that the intensity of the counting process is
bounded, we do not have time of explosion and we have a solution
defined for all time $t$ (see \cite{poisson} or \cite{QR} for all
details concerning such stochastic differential equations).
\end{pf}\\

Equations $(\ref{modified222})$ or $(\ref{modified111})$ (with
truncature) admit then a unique solution $(\rho_t)$. In Section
$3$, from convergence result, we show that these solutions are
valued in the set of states, the truncature method will be then
not necessary. Therefore solutions of $(\ref{modified222})$ and
$(\ref{modified111})$ become solutions of $(\ref{modified22})$ and
$(\ref{modified11})$ and as they are valued in the set of states
they become solutions of $(\ref{poissonseul})$ and
$(\ref{poissonseul})$.

Before to tackle the problem of convergence in Section $3$, let us
give a proposition concerning martingale problem for
$(\overline{\mathcal{A}}^J,\rho_0)$ (for some $\rho_0$) and
uniqueness of the solution for such problem. This will be namely
useful in Section $3$.

\begin{pr}\label{uniquenessmartingaleproblem}
Let $\rho_0$ be any operator. Let $\overline{\mathcal{A}}_k^J$ be
the infinitesimal generator of the process
$(\rho^{\mbox{\tiny$J$}}_t)$, solution of a truncated equation of
the form $(\ref{modified222})$ or $(\ref{modified111})$.

The process $(\rho^{\mbox{\tiny$J$}}_t)$ is then the unique
solution in distribution of the martingale problem
$(\overline{\mathcal{A}}_k^J,\rho_0)$.
\end{pr}

The fact that the solution of a stochastic differential equation
$(\ref{modified222})$ or $(\ref{modified111})$ is a solution of
the martingale problem for the corresponding infinitesimal
generator follows from Ito formula as in Proposition $3$.

This proposition means that all other solution of the martingale
problem for $(\overline{\mathcal{A}}^J,\rho_0)$ have the same
distribution of the solution $(\rho_t^{\mbox{\tiny$J$}})$ of the
associated stochastic differential equation. This result is
classical in Markov Process Generator Theory, it follows from the
pathwise uniqueness of the solutions of equations
$(\ref{modified222})$ and $(\ref{modified111})$ (see
\cite{MR838085} for a complete reference about existence and uniqueness of solutions for problems of martingale).

\section{Convergence of Discrete Quantum Trajectories}

In all this section, we consider an observable $A$ of the form
$(\ref{decompo})$ with associated subset $J$ and $I$ as in Proposition $2$. Furthermore we
consider an integer $k>1$ and the associated truncated stochastic
differential equations $(\ref{modified222})$ or
$(\ref{modified111})$.

 In this
section, we show that the discrete quantum trajectory
$(\rho_n^{\mbox{\tiny$J$}}(t))$ (describing the successive
measurements of $A$) converges in distribution to the solution of
the martingale problem for $(\overline{\mathcal{A}}_k^J,\rho_0)$
given by the solution of the corresponding truncated equations
$(\ref{modified222})$ or $(\ref{modified111})$. Next we show that
such convergence results allow to conclude that solutions of
$(\ref{modified222})$ or $(\ref{modified111})$ are valued in the
set of states.\\

Let $\rho_0$ be any initial state. In order to prove that the
discrete trajectory starting from $\rho_0$ converges in
distribution, we show at first that the finite dimensional
distributions of the discrete process
$(\rho_n^{\mbox{\tiny$J$}}(t))$ converge to the finite dimensional
distribution of the solution of the martingale problem
$(\overline{\mathcal{A}}_k^J,\rho_0)$. Secondly we show that the
discrete process $(\rho_n^{\mbox{\tiny$J$}}(t))$ is tight and the
convergence follows. For the weak convergence of finite
dimensional distributions, we use the following theorem of Ethier
and Kurtz \cite{MR838085} translated in the context of quantum
trajectories.

\begin{thm} Let $\overline{\mathcal{A}}_k^J$ be the infinitesimal
 generator of the solution of the corresponding equation $(\ref{modified222})$ or $(\ref{modified111})$. Let $(\mathcal{F}_t^n)$ be a filtration and let $(\rho_n^{\mbox{\tiny$J$}}(.))$ be
   a c\`adl\`ag $\mathcal{F}_t^n$
   adapted-process which is relatively compact (or tight). Let $\rho_0$ be any state.

Suppose that:
\begin{enumerate}
 \item{The martingale problem $(\overline{\mathcal{A}}^J_k,\rho_0)$ has a unique solution(in distribution);}
\item{$\rho^{\mbox{\tiny$J$}}_n(0)=\rho_0$;}
\item{for all $m\geq0$, for all $0\leq t_1<t_2<\ldots<t_m\leq t<t+s$, for all function
 $(\theta_i)_{i=1,\ldots, m}$ and for all $f$ in $C^2_c$ we have:\begin{equation}\label{CC}
 \lim_{n\rightarrow\infty}\mathbf{E}\left[\left(f(\rho^{\mbox{\tiny$J$}}_n(t+s))-f(\rho^{\mbox{\tiny$J$}}_n(t)
 -\int_t^{t+s}\overline{\mathcal{A}}^J_kf(\rho^{\mbox{\tiny$J$}}_n(s))ds\right)\prod_{i=1}^m\theta_i(\rho^{\mbox{\tiny$J$}}_n(t_i))\right]=0.
\end{equation}}
\end{enumerate}
Then $(\rho^{\mbox{\tiny$J$}}_n(.))$ converges in distribution to
the solution of the martingale problem for
$(\overline{\mathcal{A}}_k^J,\rho_0).$
\end{thm}

Theorem $4$ of Ethier and Kurtz imposes to have uniqueness of
solution for the martingale problem, this follows from Proposition
$4$ of Section $2$. This Theorem expresses the fact that if a subsequence of $(\rho_n(t))$ converges in distibution to a stochastic process $(Y_t)$, necessarily this process is a solution of the problem of martingale associated with $(\overline{\mathcal{A}}_k^J,\rho_0).$. Indeed, from the convergence property $(\ref{CC})$, the process $(Y_t)$ satisfies
\begin{equation}
 \mathbf{E}\left[\left(f(Y_{t+s})-f(Y_t)
 -\int_t^{t+s}\overline{\mathcal{A}}^J_kf(Y_s)\,ds\right)\prod_{i=1}^m\theta_i(\rho^{\mbox{\tiny$J$}}_n(t_i))\right]=0.
\end{equation}
As this equality is satisfied for all $m\geq0$, for all $0\leq t_1<t_2<\ldots<t_m\leq t<t+s$, for all function
 $(\theta_i)_{i=1,\ldots, m}$ and for all $f$ in $C^2_c$, this implies the martingale property of the process
$$t\rightarrow f(Y_t)-f(Y_0)-\int_0^t\mathcal{A}^Jf(Y_s)ds.$$
Hence, the uniqueness of the solution of the problem of martingale allow to conclude to the convergence of finite dimensional distributions and the tightness property allow to conclude to the convergence in distribution for stochastic processes.

 Let us deal with the application of Theorem $4$ in the context of quantum trajectories. Concerning the definition of a filtration
$(\mathcal{F}_t^n)$, we consider the natural filtration of the
discrete quantum trajectory $(\rho_n^{\mbox{\tiny$J$}}(t))$, that
is, if $r/n\leq t<(r+1)/n$ we have
$$\mathcal{F}_t^n=\sigma(\rho_n^{\mbox{\tiny$J$}}(s),s\leq t)=\sigma(\rho^{\mbox{\tiny$J$}}_p,p\leq r).$$
It is obvious that $\mathcal{F}_t^n=\mathcal{F}_{r/n}^n$.

 Assume
tightness for instance. In order to conclude, it suffices to prove
the assertion $(\ref{CC})$. As $k$ is supposed to be strictly
larger than $1$, recall that infinitesimal generators of quantum
trajectories $\mathcal{A}^J$ satisfy for all $f\in C^2_c$ and for
all states $\rho$
$$\overline{\mathcal{A}}_k^Jf(\rho)=\mathcal{A}^Jf(\rho).$$ The assertion $(\ref{CC})$ follows then from this
proposition.
\begin{pr} Let $\rho_0$ be any state. Let $(\rho_n^{\mbox{\tiny$J$}}(.))$ be
a quantum trajectory starting from $\rho_0$. Let $\mathcal{F}_t^n$
be the natural filtration of $(\rho_n^{\mbox{\tiny$J$}}(.))$. We
have:
 \begin{eqnarray}
 &&\lim_{n\rightarrow\infty}\mathbf{E}\left[\left(f(\rho^{\mbox{\tiny$J$}}_n(t+s))-f(\rho^{\mbox{\tiny$J$}}_n(t)-
 \int_t^{t+s}\overline{\mathcal{A}}_k^Jf(\rho^{\mbox{\tiny$J$}}_n(s))ds\right)\prod_{i=1}^m\theta_i(\rho^{\mbox{\tiny$J$}}_n(t_i))\right]\nonumber\\
&=&\lim_{n\rightarrow\infty}\mathbf{E}\left[\left(f(\rho^{\mbox{\tiny$J$}}_n(t+s))-f(\rho^{\mbox{\tiny$J$}}_n(t)
-\int_t^{t+s}\mathcal{A}^Jf(\rho^{\mbox{\tiny$J$}}_n(s))ds\right)\prod_{i=1}^m\theta_i(\rho^{\mbox{\tiny$J$}}_n(t_i))\right]\nonumber\\
&=&0
\end{eqnarray}
for all $m\geq0$, for all $0\leq t_1<t_2<\ldots<t_m\leq t<t+s$, for
all functions $(\theta_i)_{i=1,\ldots, m}$ and for all $f$ in $C^2_c$.
\end{pr}

\begin{pf}
 The discrete quantum trajectory $(\rho_n^{\mbox{\tiny$J$}}(t))$ is
  valued in the set of states, we then have for all $s\geq0$
$$\overline{\mathcal{A}}_k^Jf(\rho^{\mbox{\tiny$J$}}_n(s))=\mathcal{A}^Jf(\rho^{\mbox{\tiny$J$}}_n(s)).$$
Let $m\geq0$, let $0\leq t_1<t_2<\ldots<t_m\leq t<t+s$, let
$(\theta_i)_{i=1,\ldots, m}$ and let $f$ be functions in $C^2_c$, we
have
\begin{eqnarray}
 &&\mathbf{E}\left[\left(f(\rho^{\mbox{\tiny$J$}}_n(t+s))-f(\rho^{\mbox{\tiny$J$}}_n(t)-\int_t^{t+s}\mathcal{A}^Jf(\rho^{\mbox{\tiny$J$}}_n(s))ds\right)\prod_{i=1}^m\theta_i(\rho^{\mbox{\tiny$J$}}_n(t_i))\right]\nonumber\\
&=&\mathbf{E}\left[\mathbf{E}\left[\left(f(\rho^{\mbox{\tiny$J$}}_n(t+s))-f(\rho^{\mbox{\tiny$J$}}_n(t)
-\int_t^{t+s}\mathcal{A}^Jf(\rho^{\mbox{\tiny$J$}}_n(s))ds\right)\prod_{i=1}^m\theta_i(\rho^{\mbox{\tiny$J$}}_n(t_i))\right]/\mathcal{F}_t^n\right]\nonumber\\
&=&\mathbf{E}\left[\mathbf{E}\left[\left(f(\rho^{\mbox{\tiny$J$}}_n(t+s))
-f(\rho^{\mbox{\tiny$J$}}_n(t)-\int_t^{t+s}\mathcal{A}^if(\rho_n(s))ds\right)/\mathcal{F}_t^n\right]\prod_{i=1}^m\theta_i(\rho^{\mbox{\tiny$J$}}_n(t_i))\right]
\end{eqnarray}
Let $n$ be fixed, from definition of infinitesimal generators for Markov chains
(see \cite{MR2190038}) we have that
\begin{equation}\label{martingaleproperty}
 f(\rho_n^{\mbox{\tiny$J$}}(k/n))-f(\rho_0)-\sum_{j=0}^{k-1}\frac{1}{n}\mathcal{A}^J_nf(\rho^{\mbox{\tiny$J$}}_n(j/n))
\end{equation}
is a $(\mathcal{F}_{k/n}^n)$ martingale (this is the discrete
equivalent of solutions for problems of martingale for discrete processes). 

Suppose $r/n\leq t<(r+1)/n$
and $l/n\leq t+s<(l+1)/n$, we have
$\mathcal{F}_t^n=\mathcal{F}_{r/n}^n$. The random states
$\rho^{\mbox{\tiny$J$}}_n(t)$ and $\rho^{\mbox{\tiny$J$}}_n(t+s)$
satisfy then
$\rho^{\mbox{\tiny$J$}}_n(t)=\rho^{\mbox{\tiny$J$}}_n(r/n)$ and
$\rho_n^{\mbox{\tiny$J$}}(t+s)=\rho_n^{\mbox{\tiny$J$}}(l/n)$. The
martingale property $(\ref{martingaleproperty})$ implies then
\begin{eqnarray}
 &&\mathbf{E}\left[f(\rho^{\mbox{\tiny$J$}}_n(t+s))-f(\rho_n^{\mbox{\tiny$J$}}(t)\Big{/}\mathcal{F}_t^n\right]\nonumber\\
 &=&\mathbf{E}\left[f(\rho^{\mbox{\tiny$J$}}_n(l/n))-f(\rho_n^{\mbox{\tiny$J$}}(k/n)\Big{/}\mathcal{F}_{r/n}^n\right]\nonumber\\
&=&\mathbf{E}\left[\sum_{j=k}^{l-1}\frac{1}{n}\mathcal{A}^J_nf(\rho_n(j/n))\Big{/}\mathcal{F}_{r/n}^n\right]\nonumber\\
&=&\mathbf{E}\left[\int_t^{t+s}\mathcal{A}_n^Jf(\rho^{\mbox{\tiny$J$}}_n(s))ds\Big{/}\mathcal{F}_t^n\right]\nonumber\\&&
+\mathbf{E}\left[\left(t-\frac{r}{n}\right)\mathcal{A}_n^Jf(\rho^{\mbox{\tiny$J$}}_n(t))+\left(\frac{l}{n}-(t+s)\right)\mathcal{A}^J_nf(\rho^{\mbox{\tiny$J$}}_n(t+s))\Big{/}\mathcal{F}_t^n\right].
\end{eqnarray}
As a consequence, we have
\begin{eqnarray}
 &&\Bigg{\vert}\mathbf{E}\left[\left(f(\rho^{\mbox{\tiny$J$}}_n(t+s))-f(\rho^{\mbox{\tiny$J$}}_n(t)-\int_t^{t+s}
 \mathcal{A}^Jf(\rho^{\mbox{\tiny$J$}}_n(s))\right)\prod_{i=1}^m\theta_i(\rho^{\mbox{\tiny$J$}}_n(t_i))\right]
 \Bigg{\vert}\nonumber\\
&\leq&\mathbf{E}\left[\bigg{\vert}\int_t^{t+s}\Big(\mathcal{A}^J_nf(\rho_n(s))-
\mathcal{A}^Jf(\rho^{\mbox{\tiny$J$}}_n(s))\Big)ds \bigg{\vert}\right]\prod_{i=1}^m\Vert \theta_i\Vert_{\infty}\nonumber\\
&&+\mathbf{E}\left[\bigg{\vert}\left(t-\frac{[nt]}{n}\right)\mathcal{A}_nf(\rho^{\mbox{\tiny$J$}}_n(t))+
\left(\frac{[n(t+s)]}{n}-(t+s)\right)\mathcal{A}^J_nf(\rho^{\mbox{\tiny$J$}}_n(t+s)) \bigg{\vert}\right]\prod_{i=1}^m
\Vert \theta_i\Vert_{\infty}\nonumber\\
&\leq& M\sup_{\rho\in\mathcal{S}}\bigg{\vert}\mathcal{A}_n^Jf(\rho)-\mathcal{A}^Jf(\rho)\bigg{\vert}+\frac{L}{n}\sup_{\rho\in\mathcal{S}}\bigg{\vert}\mathcal{A}_n^Jf(\rho)\bigg{\vert}
\end{eqnarray}
where $M$ and $L$ are constant depending on $\Vert h_i\Vert$ and
$s$. Thanks to the condition of uniform convergence of Proposition $2$, we
obtain
\begin{equation}
 \lim_{n\rightarrow\infty}\Bigg{\vert}\mathbf{E}\left[\left(f(\rho^{\mbox{\tiny$J$}}_n(t+s))
 -f(\rho^{\mbox{\tiny$J$}}_n(t)-\int_t^{t+s}\mathcal{A}_if(\rho^{\mbox{\tiny$J$}}_n(s))\right)\prod_{i=1}^m\theta_i(\rho^{\mbox{\tiny$J$}}_n(t_i))\right]\Bigg{\vert}=0
\end{equation}
\end{pf}

Finally, in order to apply Theorem $4$ of Ethier and Kurtz, it
suffices to prove the tightness of the discrete quantum
trajectory. The following lemma will be useful when we deal with
this question.

\begin{lem} There exists a constant $K_J$ such that for all $(r,l)\in(\mathbb{N}^\star)^2$ satisfying $r<l$,
  we have almost surely:
\begin{equation}
 \mathbf{E}\left[\Vert\rho_l^{\mbox{\tiny$J$}}-\rho_r^{\mbox{\tiny$J$}}\Vert^2\Big{/}\mathcal{M}^{(n)}_r\right]
 \leq K_J\frac{l-r}{n},
\end{equation}
where
$\mathcal{M}^{(n)}_r=\mathcal{F}_{r/n}^n=\sigma\{\rho_j^{\mbox{\tiny$J$}},j\leq
r\}$.
\end{lem}

\begin{pf}
 Let us deal with the case where $J\neq\emptyset$ and $I\neq\emptyset$ (this is the most general case). Let $r<l$, we have
$$\mathbf{E}\left[\Vert\rho_l^i-\rho_r^i\Vert^2/\mathcal{M}^{(n)}_r\right]
=\mathbf{E}\left[\mathbf{E}\left[\Vert\rho_l^i-\rho_r^i\Vert^2/\mathcal{M}^{(n)}_{l-1}\right]\Big{/}\mathcal{M}^{(n)}_r
\right]$$ Hence we have
\begin{eqnarray}\label{cdf}
 \mathbf{E}\left[\Vert\rho_l^{\mbox{\tiny$J$}}-\rho_r^{\mbox{\tiny$J$}}\Vert^2/\mathcal{M}^{(n)}_{l-1}\right]
 &=&\mathbf{E}\left[\left\Vert\sum_{j=0}^p\mathcal{L}_j^{(n)}
 (\rho_{l-1}^{\mbox{\tiny$J$}})\mathbf{1}^{l+1}_j-\rho_r^{\mbox{\tiny$J$}}
 \right\Vert^2\bigg{/}\mathcal{M}^{(n)}_{l-1}\right]\nonumber\\
&=&\mathbf{E}\left[\sum_{j=0}^p\left\Vert\mathcal{L}_j^{(n)}(\rho^{\mbox{\tiny$J$}}_{l-1})
-\rho_r^{\mbox{\tiny$J$}}\right\Vert^2p^j_l(\rho^{\mbox{\tiny$J$}}_{l-1})\bigg{/}\mathcal{M}^{(n)}_{l-1}\right]\nonumber\\
&=&\mathbf{E}\left[\sum_{j=0}^p\left\Vert\mathcal{L}_j^{(n)}(\rho_{l-1}^{\mbox{\tiny$J$}})-\rho^{\mbox{\tiny$J$}}_{l-1}+
\rho^{\mbox{\tiny$J$}}_{l-1}-\rho_r^{\mbox{\tiny$J$}}\right\Vert^2p^j_l(\rho^{\mbox{\tiny$J$}}_{l-1})\bigg{/}
\mathcal{M}^{(n)}_{l-1}\right]\nonumber\\
&=&\sum_{j\in
I}\mathbf{E}\left[\left\Vert\mathcal{L}_j^{(n)}(\rho^{\mbox{\tiny$J$}}_{l-1})-\rho^{\mbox{\tiny$J$}}_{l-1}
+
\rho^{\mbox{\tiny$J$}}_{l-1}-\rho_r^{\mbox{\tiny$J$}}\right\Vert^2p^j_l(\rho^{\mbox{\tiny$J$}}_{l-1})\bigg{/}
\mathcal{M}^{(n)}_{l-1}\right]\\&&+\sum_{j\in
J\bigcup\{0\}}\mathbf{E}\left[\left\Vert\mathcal{L}_j^{(n)}(\rho^{\mbox{\tiny$J$}}_{l-1})-\rho^i_{l-1}+\rho^i_{l-1}-\rho_k^i\right\Vert^2p^j_l(\rho_{l-1})\bigg{/}\mathcal{M}^{(n)}_{l-1}\right]\nonumber\\
\end{eqnarray}
As $I$ is supposed to be not empty, for the first term of $(\ref{cdf})$ we have for all $i\in I$
$$p_{l}^i(\rho^{\mbox{\tiny$J$}}_{l-1})=\frac{1}{n}(v_i(\rho^{\mbox{\tiny$J$}}_{l-1})+\circ(1))$$
where functions $v_i(.)$ are defined in Section $2$. As
$\mathcal{L}^{(n)}_i(\rho)$ converges uniformly in
$\rho\in\mathcal{S}$, set
$$R=\sup_{j\in I}\sup_n\sup_{(\rho,\mu)\in\mathcal{S}^2}\Bigg{\{}\left\Vert\mathcal{L}_j^{(n)}(\rho)
-\mu\right\Vert^2(v_i(\rho)+\circ(1))\Bigg{\}}.$$ This constant is
finite because all the $\circ's$ are uniform in $\rho$. We have
then almost surely
$$\sum_{j\in I}\mathbf{E}\left[\left\Vert\mathcal{L}_j^{(n)}(\rho_{l-1}^{\mbox{\tiny$J$}})-\rho^{\mbox{\tiny$J$}}_{l-1}+\rho^{\mbox{\tiny$J$}}_{l-1}-\rho_r^{\mbox{\tiny$J$}}\right\Vert^2p^j_l(\rho^{\mbox{\tiny$J$}}_{l-1})\bigg{/}\mathcal{M}^{(n)}_{l-1}\right]\leq\frac{card(J)\times R}{n}$$
For the second term of $(\ref{cdf})$ we have
\begin{eqnarray}\label{cdfg}
 &&\sum_{j\in J\bigcup\{0\}}\mathbf{E}\left[\left\Vert\mathcal{L}_j^{(n)}(\rho_{l-1}^{\mbox{\tiny$J$}})
 -\rho^{\mbox{\tiny$J$}}_{l-1}+\rho^{\mbox{\tiny$J$}}_{l-1}-\rho_r^{\mbox{\tiny$J$}}
 \right\Vert^2p^j_l(\rho^{\mbox{\tiny$J$}}_{l-1})\bigg{/}\mathcal{M}^{(n)}_{l-1}\right]\nonumber\\
&=&\sum_{j\in
J\bigcup\{0\}}\mathbf{E}\left[\left\Vert\mathcal{L}_j^{(n)}(\rho_{l-1}^{\mbox{\tiny$J$}})
-\rho^{\mbox{\tiny$J$}}_{l-1}\right\Vert^2p^j_l(\rho^{\mbox{\tiny$J$}}_{l-1})\bigg{/}\mathcal{M}^{(n)}_{l-1}\right]\nonumber\\
&&+\sum_{j\in
J\bigcup\{0\}}\mathbf{E}\left[2Re\left(\langle\mathcal{L}_j^{(n)}(\rho_{l-1}^{\mbox{\tiny$J$}})-\rho_{l-1}^{\mbox{\tiny$J$}},
\rho^{\mbox{\tiny$J$}}_{l-1}-\rho_r^{\mbox{\tiny$J$}}\rangle\right)
p^j_l(\rho^{\mbox{\tiny$J$}}_{l-1})\bigg{/}\mathcal{M}^{(n)}_{l-1}\right]\nonumber\\&&+\sum_{j\in
J\bigcup\{0\}}\mathbf{E}\left[\left\Vert\rho_{l-1}^{\mbox{\tiny$J$}}-\rho^{\mbox{\tiny$J$}}_{r}
\right\Vert^2p^j_l(\rho^{\mbox{\tiny$J$}}_{l-1})\bigg{/}\mathcal{M}^{(n)}_{l-1}\right]
\end{eqnarray}
Concerning the indexes $j\in J\bigcup\{0\}$, we have
$$\mathcal{L}_j^{(n)}(\rho_{l-1}^{\mbox{\tiny$J$}})-\rho^{\mbox{\tiny$J$}}_{l-1}=\frac{1}{\sqrt{n}}(h_j(\rho^{\mbox{\tiny$J$}}_{l-1})+\circ(1)).$$
In the same way as the constant $R$, we define
$$S=\sup_{j\in J\bigcup\{0\}}\sup_n\sup_{\rho\in\mathcal{S}}\left\Vert h_j(\rho)+\circ(1)\right\Vert^2p_j(\rho).$$
For the first term of $(\ref{cdfg})$, it implies that we have
almost surely
$$\sum_{j\in J\bigcup\{0\}}\mathbf{E}\left[\left\Vert\mathcal{L}_j^{(n)}(\rho_{l-1}^{\mbox{\tiny$J$}})
-\rho^{\mbox{\tiny$J$}}_{l-1}\right\Vert^2p^j_l(\rho^{\mbox{\tiny$J$}}_{l-1})
\bigg{/}\mathcal{M}^{(n)}_{l-1}\right]\leq\frac{(card{J}+1)\times
S}{n}.$$ Furthermore, as we have $\sum_{j\in
J\bigcup\{0\}}p^j_l(\rho^{\mbox{\tiny$J$}}_{l-1})\leq 1$ almost
surely, we get almost surely
$$\sum_{j\in J\bigcup\{0\}}\mathbf{E}\left[\left\Vert\rho_{l-1}^{\mbox{\tiny$J$}}-\rho^{\mbox{\tiny$J$}}_{r}
\right\Vert^2p^j_l(\rho_{l-1})\bigg{/}\mathcal{M}^{(n)}_{l-1}\right]\leq\mathbf{E}\left[\left\Vert\rho_{l-1}^{\mbox{\tiny$J$}}
-\rho^{\mbox{\tiny$J$}}_{r}\right\Vert^2\bigg{/}\mathcal{M}^{(n)}_{l-1}\right].$$

 For the second term of $(\ref{cdfg})$, we have
\begin{eqnarray*}
 &&\sum_{j\in J\bigcup\{0\}}\mathbf{E}\left[2Re\left(\left\langle\mathcal{L}_j^{(n)}(\rho_{l-1}^{\mbox{\tiny$J$}})
 -\rho_{l-1}^{\mbox{\tiny$J$}},\rho^{\mbox{\tiny$J$}}_{l-1}-\rho_r^{\mbox{\tiny$J$}}\right\rangle\right)
 p^j_l(\rho^{\mbox{\tiny$J$}}_{l-1})\bigg{/}\mathcal{M}^{(n)}_{l-1}\right]\\
&=&\mathbf{E}\left[2Re\left(\left\langle\sum_{j\in
J\bigcup\{0\}}\left(\mathcal{L}_j^{(n)}(\rho_{l-1}^{\mbox{\tiny$J$}})-\rho_{l-1}^{\mbox{\tiny$J$}}
\right)p^j_l(\rho^{\mbox{\tiny$J$}}_{l-1})\,,\,\rho^{\mbox{\tiny$J$}}_{l-1}-\rho_r^{\mbox{\tiny$J$}}\right\rangle\right)
\bigg{/}\mathcal{M}^{(n)}_{l-1}\right]
\end{eqnarray*}
Let us treat this term. As in the proof of Proposition $2$ concerning  infinitesimal generators, with the
asymptotic of $\mathcal{L}_i^{(n)}$ in this situation, we have
uniformly in $\rho\in\mathcal{S}$:
$$\sum_{j\in J\bigcup\{0\}}\left(\mathcal{L}_j^{(n)}(\rho)-\rho\right)p^j_l(\rho)=\frac{1}{n}(H(\rho)+\circ(1)),$$
since the terms in $1/\sqrt{n}$ disappear by summing over $j\in
J\bigcup\{0\}$. As a consequence, by defining the finite constant
$W$ as
$$W=\sup_n\sup_{(\rho,\mu)\in\mathcal{S}^2}\Bigg{\{}\left\vert2Re\left(\left\langle\sum_{j\in J\bigcup\{0\}}n
\left(\mathcal{L}_j^{(n)}(\rho)-\rho\right)p^j_l(\rho)\,,\,\rho-\mu\right\rangle\right)
\right\vert\Bigg{\}},$$ we have then almost surely
$$\sum_{j\in J\bigcup\{0\}}\mathbf{E}\left[2Re\left(\left\langle\mathcal{L}_j^{(n)}(\rho_{l-1}^{\mbox{\tiny$J$}})-\rho_{l-1}^{\mbox{\tiny$J$}}
,\rho^{\mbox{\tiny$J$}}_{l-1}-\rho_r^{\mbox{\tiny$J$}}\right\rangle\right)
p^j_l(\rho^{\mbox{\tiny$J$}}_{l-1})\bigg{/}\mathcal{M}^{(n)}_{l-1}\right]\leq
\frac{W}{n}$$
Let us stress that the constant $W$ are independent of $l$ and $r$. Therefore, we can conclude that  there exists a constant $K_J$ such that for all
$r<l$, we have almost surely
\begin{eqnarray}
 \mathbf{E}\left[\left\Vert\rho^{\mbox{\tiny$J$}}_l-\rho^{\mbox{\tiny$J$}}_r\right\Vert^2/\mathcal{M}^{(n)}_{l-1}\right]
 \leq \frac{K_J}{n}+\mathbf{E}\left[\left\Vert\rho^{\mbox{\tiny$J$}}_{l-1}-\rho^{\mbox{\tiny$J$}}_r\right
 \Vert^2/\mathcal{M}^{(n)}_{l-1}\right].
\end{eqnarray}
It implies that almost surely
\begin{equation}
\mathbf{E}\left[\left\Vert\rho^{\mbox{\tiny$J$}}_l-\rho_r^{\mbox{\tiny$J$}}
\right\Vert^2/\mathcal{M}^{(n)}_{k}\right]\leq\frac{K_J}{n}+\mathbf{E}\left[\left\Vert\rho^{\mbox{\tiny$J$}}_{l-1}-
\rho^{\mbox{\tiny$J$}}_r\right\Vert^2/\mathcal{M}^{(n)}_{r}\right].
\end{equation}
Thus by conditioning successively by $\mathcal{M}^{(n)}_{l-i}$
with $i\in\{2,\ldots,l-r\}$ and by induction, we can show
$$\mathbf{E}\left[\left\Vert\rho^{\mbox{\tiny$J$}}_l-\rho^{\mbox{\tiny$J$}}_r\right\Vert^2/\mathcal{M}^{(n)}_{r}\right]
\leq\frac{K_J(l-r)}{n}.$$ The same results hold when $J=\emptyset$
or $I=\emptyset$ by similar computations.
\end{pf}\\

This lemma implies the following proposition which expresses the
tightness property of the discrete quantum trajectory.

\begin{pr} Let $(\rho^{\mbox{\tiny$J$}}_n(t))$ be the quantum trajectory. There
exists some constant $Z_J$ such that for all $t_1<t<t_2$:
\begin{equation}\label{tight}
\mathbf{E}\left[\Vert\rho_n(t_2)-\rho_n(t)\Vert^2\Vert\rho_n(t)-\rho_n(t_1)\Vert^2\right]\leq Z_J(t_2-t_1)^2.
\end{equation}
Therefore, the discrete quantum trajectory
$(\rho^{\mbox{\tiny$J$}}_n(t))$ is tight.
\end{pr}

\begin{pf}
 The inequality $(\ref{tight})$ implies the tightness of $(\rho_n(t))$ (see \cite{MR1700749}). Let us prove $(\ref{tight})$. It is worth noticing that $\mathcal{M}^{(n)}_k=\mathcal{F}_{k/n}^n$ where $\mathcal{F}_t^n$ is the natural filtration of $(\rho_n^{\mbox{\tiny$J$}})$. Thanks to the previous lemma we then have:
\begin{eqnarray*}
 &&\mathbf{E}\left[\Vert\rho^{\mbox{\tiny$J$}}_n(t_2)-\rho^{\mbox{\tiny$J$}}_n(t)\Vert^2\Vert\rho^{\mbox{\tiny$J$}}_n(t)-\rho_n(t_1)\Vert^2\right]\\&=&
\mathbf{E}\left[\Vert\rho^{\mbox{\tiny$J$}}_n([nt_2])-\rho^{\mbox{\tiny$J$}}_n([nt])\Vert^2\Vert\rho^{\mbox{\tiny$J$}}_n([nt])-\rho^{\mbox{\tiny$J$}}_n([nt_1])\Vert^2\right]\\&=&\mathbf{E}\left[\mathbf{E}\left[\Vert\rho^{\mbox{\tiny$J$}}_n([nt_2])-\rho^{\mbox{\tiny$J$}}_n([nt])\Vert^2\Vert\rho^{\mbox{\tiny$J$}}_n([nt])-\rho^{\mbox{\tiny$J$}}_n([nt_1])\Vert^2/\mathcal{F}_{[nt]/n}^n\right]\right]\\
&=&\mathbf{E}\left[\mathbf{E}\left[\Vert\rho^{\mbox{\tiny$J$}}_n([nt_2])-\rho^{\mbox{\tiny$J$}}_n([nt])\Vert^2/\mathcal{F}_{[nt]/n}^n\right]\Vert\rho^{\mbox{\tiny$J$}}_n([nt])-\rho_n([nt_1])\Vert^2\right]\\
&\leq&\mathbf{E}\left[\frac{K_J([nt_2]-[nt])}{n}\Vert\rho^{\mbox{\tiny$J$}}_n([nt])-\rho^{\mbox{\tiny$J$}}_n([nt_1])\Vert^2\right]\\
&\leq&\frac{K_J([nt_2]-[nt])}{n}\mathbf{E}\left[\mathbf{E}\left[\Vert\rho^{\mbox{\tiny$J$}}_n([nt])-\rho^{\mbox{\tiny$J$}}_n([nt_1])\Vert^2/\mathcal{F}_{[nt_1]/n}^n\right]\right]\\
&\leq&\frac{K_J([nt_2]-[nt])}{n}\frac{K_J([nt]-[nt_1])}{n}\\
&\leq&Z_J(t_2-t_1)^2,
\end{eqnarray*}
with $Z_J=4(K_J)^2$ and the result follows.
\end{pf}\\

Therefore we have proved the tightness property of discrete quantum
trajectories and we can now express the final theorem.

\begin{thm}Let $A$ be an observable on $\mathcal{H}=\mathbb{C}^{N+1}$ with spectral decomposition \begin{equation}A=\sum_{i=0}^p\lambda_iP_i=\sum_{i\in I}\lambda_iP_i+\sum_{j\in
J\bigcup0}\lambda_jP_j,\end{equation} where: \begin{enumerate}
\item{For $i\in\{0,\ldots, p\}$ the operators
$P_i=(p^i_{kl})_{0\leq k,l\leq N}$ are the eigen-projectors of $A$
(satisfying $p^0_{00}\neq0)$} \item{ The sets $I$ and $J$ satisfy that
$I=\{i\in\{1,\ldots,p\}/p^i_{00}=0\}$ and
$J=\{1,\ldots,p\}\setminus I$.} \end{enumerate} Let $\rho_0$ be a
state on $\mathcal{H}_0$. Let $(\rho_n^{\mbox{\tiny$J$}}(t))$ be
the discrete quantum trajectory describing the repeated quantum
measurements of $A$ and starting with $\rho_0$ as initial state.
\begin{enumerate}
\item{Suppose $J=\emptyset$. Then the discrete quantum trajectory
$(\rho_n^{\mbox{\tiny$J$}}(t))$ converges in distribution in
$\mathcal{D}[0,T)$ for all $T$ to the solution
$(\rho_t^{\mbox{\tiny$J$}})$ of the stochastic differential
equation $(\ref{modified222})$. Therefore the process
$(\rho_t^{\mbox{\tiny$J$}})$ takes values in the set of states on
$\mathcal{H}_0$. The discrete quantum trajectory $(\rho_n^{\mbox{\tiny$J$}}(t))$ converges then in
distribution to the unique solution of the following
jump-diffusion Belavkin equation
\begin{eqnarray}\label{AAa}
\rho^{\mbox{\tiny$J$}}_t=\rho_0+\int L(\rho^{\mbox{\tiny$J$}}_{s-})ds+\sum_{i=1}^p\int_0^t\int_{\mathbb{R}}g_i(\rho^{\mbox{\tiny$J$}}_{s-})\mathbf{1}_{0<x<v_i(\rho^{\mbox{\tiny$J$}}_{s-})}\left[N_i(dx,ds)-dxds\right]
\end{eqnarray}}
\item{Suppose $J\neq\emptyset$. Then the discrete quantum trajectory
$(\rho_n^{\mbox{\tiny$J$}}(t))$ converges in distribution in
$\mathcal{D}[0,T)$ for all $T$ to the solution
$(\rho_t^{\mbox{\tiny$J$}})$ of the stochastic differential
equation $(\ref{modified111})$. The process
$(\rho_t^{\mbox{\tiny$J$}})$ takes values in the set of states on
$\mathcal{H}_0$. The discrete quantum trajectory $(\rho_n^{\mbox{\tiny$J$}}(t))$ converges then in
distribution to the unique solution of the following
jump-diffusion Belavkin equation
\begin{eqnarray}\label{BBb}
 \rho^{\mbox{\tiny$J$}}_t&=&\rho_0+\int L(\rho^{\mbox{\tiny$J$}}_{s-})ds+\sum_{i\in J\bigcup\{0\}}h_i(\rho^{\mbox{\tiny$J$}}_{s-})dW_i(s)\nonumber\\&&+\sum_{i\in I}\int_0^t\int_{\mathbb{R}}g_i(\rho^{\mbox{\tiny$J$}}_{s-})\mathbf{1}_{0<x<v_i(\rho^{\mbox{\tiny$J$}}_{s-})}\left[N_i(dx,ds)-dxds\right].
\end{eqnarray}}
\end{enumerate}

Furthermore the processes defined by
\begin{equation}
 \tilde{N}^i_t=\int_0^t\int_{\mathbb{R}}\mathbf{1}_{0<x<v_i(\rho^{\mbox{\tiny$J$}}_{s-})}N_i(dx,ds)
\end{equation}
are counting processes with stochastic intensities
$$t\rightarrow\int_0^tv_i(\rho^{\mbox{\tiny$J$}}_{s-})ds.$$
\end{thm}

As in Theorem $3$, the last assertion concerning the counting
processes of Theorem $5$ follows from properties of Poisson Point
processes $N_i$. It means actually that processes defined by
\begin{equation}
\int_0^t\int_{\mathbb{R}}\mathbf{1}_{0<x<v_i(\rho^{\mbox{\tiny$J$}}_{s-})}N_i(dx,ds)-\int_0^tv_i(\rho^{\mbox{\tiny$J$}}_{s-})ds
\end{equation}
are martingale with respect to the natural filtration of
$(\rho^{\mbox{\tiny$J$}}_{t})$ (see \cite{},\cite{}). Let us prove
the convergence results of Theorem $5$.

\begin{pf} In all cases, the convergence result follows from Theorem $4$ and Propositions $5$ for the finite dimensional distribution convergence and from proposition $6$ for the tightness. In order to finish the proof of this theorem, we have to prove that solutions of stochastic differential equations $(\ref{modified222})$ and $(\ref{modified111})$ takes values in the set of states. It is given by the convergence in distribution.

Indeed, let $(\rho_n^{\mbox{\tiny$J$}}(t))$ be converging to the
corresponding solution $(\rho_t^{\mbox{\tiny$J$}})$ of equation
$(\ref{modified222})$ or $(\ref{modified111})$, we have to prove
that this solution is self-adjoint, positive with trace $1$. By
using the convergence in distribution, we have for all
$z\in\mathbb{C}^2$
\begin{eqnarray}
 &&\rho^{\mbox{\tiny$J$}}_n(t)-(\rho^{\mbox{\tiny$J$}}_n(t))^\star\,\,\mathop{\longmapsto}^{\mathcal{D}}_{n\rightarrow\infty}\,\,\rho^{\mbox{\tiny$J$}}_t-(\rho^{\mbox{\tiny$J$}}_t)^\star\\
&&\hspace{1,1cm}Tr\left[\rho^{\mbox{\tiny$J$}}_n(t)\right]\,\,\mathop{\longmapsto}^{\mathcal{D}}_{n\rightarrow\infty}\,\,
Tr\left[\rho^{\mbox{\tiny$J$}}_t\right]\\
&&\hspace{0,9cm}\langle z,\rho_n^{\mbox{\tiny$J$}}(t)z\rangle\,\,\mathop{\longmapsto}^{\mathcal{D}}_{n\rightarrow\infty}\,\,\langle z,\rho_t^{\mbox{\tiny$J$}}z\rangle
\end{eqnarray}
where $\mathcal{D}$ denotes the convergence in distribution for processes. As $(\rho_n^{\mbox{\tiny$J$}}(t))$ takes values in the set of states, we have almost surely for all $t$ and all $z\in\mathbb{C}^2$
$$\rho^{\mbox{\tiny$J$}}_n(t)-(\rho^{\mbox{\tiny$J$}}_n(t))^\star=0,\,\,\,\,Tr\left[\rho^{\mbox{\tiny$J$}}_t\right]=1,\,\,\,\,\langle z,\rho_n^{\mbox{\tiny$J$}}(t)z\rangle\geq0.$$
These properties are conserved at the limit in distribution and the process $(\rho_t^{\mbox{\tiny$J$}})$ takes then values in the set of states. The proof of Theorem $5$ is then complete.
\end{pf}\\

This theorem is then a mathematical and physical justification of
stochastic models of continuous time quantum measurement theory.
Let us stress that in general it is difficult to prove that
equations $(\label{AA})$ and $(\label{BB})$ admit a unique
solution which takes values in the set of states. One can notice
that such equations preserve the self adjoint and trace
properties. Concerning the positivity property, it is far from
being obvious and it points out the importance of the convergence
result.

Let us conclude this article with some remarks concerning these
continuous stochastic models.

The first remark concerns the average of solutions of
$(\ref{AAa})$ or $(\ref{BBb})$. Let $(\rho_t^{\mbox{\tiny$J$}})$
be a solution of $(\ref{AAa})$ or $(\ref{BBb})$. In all cases, we
have
\begin{equation}\label{master}
 \mathbf{E}[\rho^{\mbox{\tiny$J$}}_t]=\int_0^tL(\mathbf{E}[\rho^{\mbox{\tiny$J$}}_s])ds.
\end{equation}
This follows namely from martingale property of the Brownian
motion and counting processes. This remark concerning
$(\ref{master})$ means that the function
$$t\rightarrow\mathbf{E}[\rho^{\mbox{\tiny$J$}}_t]\,,$$
is the solution of the ordinary differential equation
$$d\mu_t=L(\mu_t)dt.$$
This equation is called the ``Master Equation'' in quantum mechanics
and describes the evolution of the reference state of the small
system $\mathcal{H}_0$ without measurement. In average continuous
quantum trajectories evolve then as the solution of the Master
equation (for all measurement experiences).

The second remark concerns the classical Belavkin equations
$(\ref{diffff})$ and $(\ref{jump-equation})$. In \cite{} and
\cite{}, it is shown that such continuous model are justified from
convergence of stochastic integral and random coupling method (it
does not use infinitesimal generators theory) . With Theorem $5$,
we recover these equations by considering the case where the
measured observable $A$ is of the form
$A=\lambda_0P_0+\lambda_1P_1$. Indeed in this case, we just have
one noise at the limit as in the classical case.

The last remark concerns the uniqueness of a solution of the
martingale problems. In this article, we have made an
identification with the set of operators on $\mathcal{H}_0$ and
$\mathbb{R}^P$ in order to introduce definition of infinitesimal
generators and notion of martingale problem (see Section $2$,
Definition $2$). As observed, the infinitesimal generators of
quantum trajectories can be written in term of the partial
derivative in the following way
\begin{eqnarray}\frac{1}{2}\sum_{i\in J\bigcup\{0\}}D^2_\rho f(h_i(\rho),h_i(\rho))&=&\frac{1}{2}\sum_{i,j=1}^Pa_{ij}(\rho)\frac{\partial f(\rho)}{\partial \rho_i\partial \rho_j}\\
 D_\rho f(L(\rho)&=&\sum_{i=1}^Pb_i(\rho)\frac{\partial f(\rho)}{\partial \rho_i}
\end{eqnarray}
by expanding the differential terms. The matrix $a(.)=(a_{ij}(.)$
is a semi definite matrix. Let $W$ be a $P$ dimensional Brownian
motion, the solution of the problem of martingale can also be
expressed as the solution of
\begin{eqnarray}\label{classic}
 \rho^{\mbox{\tiny$J$}}_t&=&\rho_0+\int_0^tL(\rho^{\mbox{\tiny$J$}}_{s-})ds+\int_0^t\sigma(\rho^{\mbox{\tiny$J$}}_s)dW_s\nonumber\\&&
+\sum_{k\in I}\int_0^t\int_{\mathbb{R}}
\tilde{g}_k(\rho^{\mbox{\tiny$J$}}_{s-})
\mathbf{1}_{0<x<Re(v_k(\rho^{\mbox{\tiny$J$}}_{s-}))}\left[N_k(dx,ds)-dxds\right],
\end{eqnarray}
where $\sigma(.)$ is as matrix defined by
$\sigma(.)\sigma^t(.)=a(.)$ (see \cite{}). Let us stress that, in
this description we deal with a $P$ dimensional Brownian motion
corresponding to the dimension of $\mathbb{R}^P$ (which depends
only on the dimension of $\mathcal{H}_0$) whereas in Theorem $5$
we consider a $(p+1)$-dimensional Brownian motion corresponding to
the number of eigenvalues (which only depends on the dimension of
the interacting quantum system $\mathcal{H}$). As a consequence
from uniqueness of martingale problem (Proposition $4$) we have
two different descriptions of continuous quantum trajectories, but
they are the same as regards their distributions.
\nocite{*}
\bibliography{bibliojumpdiffusion}
\bibliographystyle{plain}

\end{document}